\documentclass[12pt]{article}
\usepackage{graphicx}
\usepackage{amsmath,amsthm,amssymb,enumerate}
\usepackage{euscript,mathrsfs}
\usepackage{url}
\usepackage[left=2cm,right=2cm,top=3.5cm,bottom=3.5cm]{geometry}
\usepackage{color}

\usepackage{soul}

\catcode`\@=11 \@addtoreset{equation}{section}

\catcode`\@=12

\allowdisplaybreaks

\newcommand{\tc}{\color{red}}

\newtheorem{Theorem}{Theorem}[section]
\newtheorem{Proposition}[Theorem]{Proposition}
\newtheorem{Lemma}[Theorem]{Lemma}
\newtheorem{Corollary}[Theorem]{Corollary}

\theoremstyle{definition}
\newtheorem{Definition}[Theorem]{Definition}

\newtheorem{Remark}[Theorem]{Remark}

\newcommand{\bTheorem}[1]{
\begin{Theorem} \label{T#1} }
\newcommand{\eT}{\end{Theorem}}

\newcommand{\bProposition}[1]{
\begin{Proposition} \label{P#1}}
\newcommand{\eP}{\end{Proposition}}

\newcommand{\bLemma}[1]{
\begin{Lemma} \label{L#1} }
\newcommand{\eL}{\end{Lemma}}

\newcommand{\bCorollary}[1]{
\begin{Corollary} \label{C#1} }
\newcommand{\eC}{\end{Corollary}}

\newcommand{\bRemark}[1]{
\begin{Remark} \label{R#1} }
\newcommand{\eR}{\end{Remark}}

\newcommand{\bDefinition}[1]{
\begin{Definition} \label{D#1} }
\newcommand{\eD}{\end{Definition}}

\newcommand{\Del}{\Delta_x}

\newcommand{\bfphi}{\boldsymbol{\varphi}}

\newcommand{\bFormula}[1]{
\begin{equation} \label{#1}}
\newcommand{\eF}{\end{equation}}

\newcommand{\Ov}[1]{\overline{#1}}

\newcommand{\DC}{C^\infty_c}

\newcommand{\aleq}{\stackrel{<}{\sim}}

\newcommand{\vr}{\varrho}

\newcommand{\vue}{\vu_\ep}

\newcommand{\vu}{\vc{u}}
\newcommand{\vm}{\vc{m}}

\newcommand{\vc}[1]{{\bf #1}}

\newcommand{\Div}{{\rm div}_x}
\newcommand{\Grad}{\nabla_x}

\newcommand{\dx}{\,{\rm d} {x}}

\newcommand{\dt}{\,{\rm d} t }

\newcommand{\intO}[1]{\int_{\Omega} #1 \ \dx}

\newcommand{\vv}{\vc{v}}

\newcommand{\D}{{\rm d}}

\newcommand{\ep}{\varepsilon}

\newcommand{\R}{\mathbb{R}}

\def\softd{{\leavevmode\setbox1=\hbox{d}%
          \hbox to 1.05\wd1{d\kern-0.4ex{\char039}\hss}}}
\definecolor{Cgrey}{rgb}{0.85,0.85,0.85}
\definecolor{Cblue}{rgb}{0.50,0.85,0.85}
\definecolor{Cred}{rgb}{1,0,0}
\definecolor{fancy}{rgb}{0.10,0.85,0.10}

\newcommand\Cbox[2]{%
    \newbox\contentbox%
    \newbox\bkgdbox%
    \setbox\contentbox\hbox to \hsize{%
        \vtop{
            \kern\columnsep
            \hbox to \hsize{%
                \kern\columnsep%
                \advance\hsize by -2\columnsep%
                \setlength{\textwidth}{\hsize}%
                \vbox{
                    \parskip=\baselineskip
                    \parindent=0bp
                    #2
                }%
                \kern\columnsep%
            }%
            \kern\columnsep%
        }%
    }%
    \setbox\bkgdbox\vbox{
        \color{#1}
        \hrule width  \wd\contentbox %
               height \ht\contentbox %
               depth  \dp\contentbox
        \color{black}
    }%
    \wd\bkgdbox=0bp%
    \vbox{\hbox to \hsize{\box\bkgdbox\box\contentbox}}%
    \vskip\baselineskip%
}

\date{}

\begin{document}


\title{Existence and stability of dissipative turbulent solutions to a simple bi-fluid model of compressible fluids}

\author{Bumja Jin
\thanks{ {The work of the first author was partially supported by { NRF-2019R1A2C1086070.}}}   \and
Young-Sam Kwon
\thanks{ {The work of the second author was partially supported by  NRF2020R1F1A1A01049805.}}
\and
\v S\' arka Ne\v casov\' a \thanks{ \v S. N.  has been supported by the Czech Science Foundation
(GA\v CR) project GA19-04243S. The Institute of Mathematics, CAS is
supported by RVO:67985840. }
\and
Anton\'{i}n Novotn{\' y}
\thanks{The work of the fourth author was partially supported by the distinguished Edurad \v Cech visiting program at the
Institute of Mathematics of the Academy of Sciences of  the Czech Republic.}}

\maketitle

\centerline{Department of Mathematics Education, Mokpo National University,
}

\centerline{Muan 534-729, South Korea}
\medskip

\centerline{Department of Mathematics, Dong-A University}
\centerline{Busan 49315, Republic of Korea, ykwon@dau.ac.kr}

\medskip
\centerline{Institute of Mathematics of the Academy of Sciences of the Czech Republic,}
\centerline{\v Zitn\'a 25, 115 67 Czech Republic, matus@math.cas.cz}

\medskip
\centerline{University of Toulon, IMATH, EA 2134,  BP 20139}
\centerline{839 57 La Garde, France, novotny@univ-tln.fr}

\begin{abstract}
Following Abbatiello et al. [ DCCDS-A (41), 2020], we introduce
dissipative turbulent solutions to a simple model of a mixture of two non interacting compressible fluids { filling
a bounded domain with general
non zero inflow/outflow boundary conditions.}
We prove existence of such solutions
for all adiabatic coefficients $\gamma>1$,  their
compatibility with classical solutions, the relative energy inequality, and the weak strong uniqueness principle in this class. The class of dissipative turbulent solutions is so far the largest class of generalized solutions which still enjoys the weak strong uniqueness property.
\end{abstract}


\bigskip

\tableofcontents

\section{Introduction}
\label{P}

The most simple system suggested as a ``toy'' problem to get a better insight into the complex mathematics in the multi-fluid modeling of compressible fluids is the following bi-fluid model for scalar
density fields $R=R(t,x)\ge 0$, $Z=Z(t,x)\ge 0$,  and vector common velocity
field $\vu=\vu(t,x)\in\R^d$ ($t\in I$, $I=(0,T)$,
$x\in\Omega$)  consisting of
\begin{enumerate}
\item {\it Conservation of mass for the species}
\begin{equation} \label{P1}
\partial_t R + \Div (R \vu) = 0; \; \partial_t Z + \Div (Z \vu) = 0\;\mbox{in $Q_T=I\times\Omega$};
\end{equation}
\item {\it Balance of linear momentum}
\begin{equation} \label{P2}
\partial_t ((R+Z) \vu) + \Div ((R+Z) \vu \otimes \vu) + \Grad P(R,Z)) = \Div \mathbb{S}\; \mbox{in $Q_T=I\times\Omega$};
\end{equation}
\item {\it Balance of energy:}
\begin{equation} \label{P3}
\partial_t \left( \frac{1}{2} (R+Z) |\vu|^2 + { H(R,Z)}\right) +
\Div \left[ \left( \frac{1}{2} (R+Z) |\vu|^2 + H(R,Z)  \right) \vu \right]
\end{equation}
$$
= \Div \left( \mathbb{S} \cdot \vu \right)
- \mathbb{S} : \Grad \vu\; \mbox{in $Q_T=I\times\Omega$},
$$
where $P(R,Z)$ is the so--called  Helmholtz function (pressure potential) related to the pressure $P=P(R,Z)$,
\begin{equation} \label{P4}
H(R,Z)=R\int_1^R\frac{P(s, s\frac Z R)}{s^2}{\rm d}s,\;\mbox{if $R>0$},\; H(0,Z)=0.
\end{equation}

\item
We suppose that the fluid  is contained in a bounded Lipschitz domain $\Omega \subset \R^d$, $d = 2,3$, with general {\it inflow--outflow boundary conditions}\footnote{We suppose without loss of generality that the boundary {  data} are restrictions to $\partial\Omega$ of functions defined on $\overline\Omega$.}

\begin{equation} \label{P5}
\vu|_{\partial \Omega} = \vu_B|_{\partial\Omega},\ R|_{\Gamma^{\rm in}} = R_B|_{\Gamma^{\rm in}},\ Z|_{\Gamma^{\rm in}} = Z_B|_{\Gamma^{\rm in}},\
\Gamma^{\rm in} = \left\{ x \in \partial \Omega \ \Big| \ \vu_B \cdot \vc{n} < 0 \right\},
\end{equation}
where $\vc{n}$ is the outer normal vector to $\partial \Omega$. It is to be noticed that no boundary { densities} are prescribed at
$$
\Gamma^{\rm out}=\partial\Omega\setminus \Gamma^{\rm in}
$$
in agreement with the nature of the equations (\ref{P1}).

For the sake of simplicity, we consider Newtonian fluids, i.e.
\begin{equation}\label{P6}
 \mathbb{S}=\mathbb{S}(\Grad\vu): = \mu( {\Grad \vu + \Grad \vu^T}) + \lambda \Div \vu \mathbb{I},\;\mu>0,\;
 \lambda+\frac 2 d\mu>0.
\end{equation}

\item Finally, we add to the system initial conditions
\begin{equation}\label{P7}
 R(0)=R_0,\; Z(0)=Z_0,\; (R+Z)\vu (0):=\vc m_0=(R_0+Z_0)\vu_0
\end{equation}
\end{enumerate}

The goal of { this paper  is to define} as weak as possible
solution
to system (\ref{P1}--\ref{P6}) (which we will call {\it dissipative turbulent solution}), which still {enjoys } the following three fundamental
properties:
\begin{enumerate}
\item {\it Existence:} Dissipative turbulent solutions exist on an arbitrary large time interval for any finite energy initial data.
\item {\it Compatibility:} If the dissipative turbulent solution is sufficiently continuously differentiable then it is a classical solution.
\item {\it Weak--strong uniqueness:} Any dissipative turbulent solution
coincides with the strong solution of the same problem emanating from the same initial and boundary data as long as the latter exist.
\end{enumerate}

It is well known from the mono-fluid theory that an object of this
type is very convenient for many applications ranging from rigorous investigation of singular limits through { dimension} reduction to
investigation of convergence and error estimates for numerical schemes
to problem (\ref{P1}--\ref{P6}).

In addition to the above mentioned favorable features of these solutions,
they also have a perfect physical interpretation. In the generalized weak formulation {we let appear} a positive semi-definite tensor $\mathfrak{R}$
(cf. Definition \ref{DW1} later) which can be interpreted
as Reynolds stress in the {\it turbulence and acoustic modeling}. Indeed, this is the same tensor which appears
as the source term in the Lignthill's acoustic analogy, cf. \cite{LH1}, \cite{LH2}.

To perform this program we shall rely on the concept of {\it dissipative solutions} with Reynolds defect introduced within the context of mono-fluid theory in \cite{AbbFei2} (incompressible fluids) and
\cite{AbFeNo} (compressible fluids). In the mono-fluid theory, this
concept is known as the most weak concept of solutions enjoying the
three { fundamental} properties stated above.

We finish this  { introductory} section with some bibliographic remarks. Existence of weak solutions to system (\ref{P1})--(\ref{P4}) with homogenous boundary conditions has been obtained by Vasseur et al. \cite{VWY} (revisited later in \cite{Wen}). The same problem with general non homogenous boundary data is investigated by Kwon et al. \cite{KwNeNo}.
The compactness argument discovered in \cite{VWY} has been generalized in
\cite{AN-MP} by using the philosophy of \cite{3MNPZ}. This argument
combined with Lions' compactness argument \cite{LI4}, opened the way to treat more realistic multi-fluid compressible models, see \cite{AN-MP}, \cite{NoSCM}, \cite{KwNeNo} with algebraic or differential { closure}. Prior
to this results, only 1-d equations or quasi-stationary multi fluid systems could be treated: a good sample of such  studies is the paper by Evje \cite{EV2} ($1-d$) and by Bresch et al. \cite{BRZA}. An interesting overview of these type of models from the point of view of mathematical physics  is the review paper by Bresch et al. \cite{BreschMF}.

The results on weak strong uniqueness  have their sources in the relative energy method \cite{Dafermos} adapted to viscous compressible fluids in \cite{FeJiNo} (homogenous boundary conditions) and in \cite{KwNo} (general in/out-flow boundary conditions). { Since the paper by Gwiazda et al. \cite{FeGw}}, it is known, that the weak strong uniqueness principle holds in larger classes than in the class
of weak solutions. Such observation has not only an academic impact: it has a practical
impact e.g. on the investigation of convergence and  error estimates for numerical schemes. In this respect larger class means less conditions
on the structure of the numerical scheme and allows more applications.
The class of solutions we considered in this paper is so far the largest one,
where still the weak strong uniqueness holds in the case of compressible Navier-Stokes equations, see \cite{AbFeNo}. This paper shows weak-strong uniqueness for the bi-fluid model (\ref{P1})-\eqref{P7} in this class.

{ The paper is organized as follows. In the next section, we define {\it (weak) dissipative (turbulent) solutions} and state the main results: Teorem \ref{ET1} (existence), Theorem \ref{ET2}
(compatibility), Theorem \ref{ET3} (relative energy inequality),
Theorem \ref{ET4} (weak-strong uniqueness). The following Sections
3-6 are devoted to the proofs of these theorems. Finally in Appendix, we recall some specific tools needed in the proofs, for reader's convenience. }

\section{Main results}\label{M}

{
Throughout the paper, we use the standard notation for Lebesgue, Sobolev and Bochner spaces, see e.g. the book of Evans \cite{Ev}.
Further, we denote by $C_{\rm weak}(\overline I;X)$ ($X$ a Banach space) the vector subspace of $L^\infty(I;X)$ of functions $f$
defined everywhere on $\overline I$ such that
for all $\eta\in X^*$, $<\eta,f>_{X^*,X}\in C(\overline I)$.
The symbol ${\cal M}(\overline\Omega)$ denotes the set of signed Radon measures on $\overline\Omega$. Symbol
$\mathcal{M}^+(\Ov{\Omega}; \R^{d \times d}_{\rm sym})$ denotes the set of all positively semi--definite tensor valued Radon measures $\mathfrak{R}$ on $\Ov{\Omega}$. This means
$\mathfrak{R}=(\mathfrak{R}_{ij})_{i,j=1,\ldots,d}$, where:
1) $\mathfrak{R}_{ij}$ is a signed Radon measure on $\overline\Omega$; 2)$\mathfrak{R}_{ij}=\mathfrak{R}_{ji}$; 3)
for all $0\neq \xi$, $\xi^T\mathfrak{R}\xi$ is poitive
Radon measure on $\Ov{\Omega}$. Finaly the Bochner type spaces
$L^\infty_{weak-*}(I;X)$ are defined in Appendix, see Lemma \ref{Lemma3}.
}

\subsection{Definition of dissipative turbulent solutions}\label{W}

Motivated by \cite[Section 2]{AbFeNo}, we introduce the dissipative
turbulent solutions to problem \eqref{P1}--\eqref{P7} as follows.

\begin{Definition}\label{DW1}
We say that a triplet $(R,Z,\vu)$,  belonging to the class
$$
{ R,Z \in C_{\rm weak}(\overline I; L^\gamma(\Omega))} \cap L^\gamma(0,T; L^\gamma(\partial \Omega; |\vu_B \cdot \vc{n}|
{\rm d}S_x))\,\mbox{with some $\gamma>1$},
$$
\begin{equation}\label{Classw}
R\ge 0, \,Z\ge 0,\,\vu-\vu_B\in L^2(I, W^{1,2}_0(\Omega;\R^d)),\;
(R+Z) \vu \in{ C_{\rm weak}([0,T]; L^{\frac{2 \gamma}{\gamma + 1}}(\Omega; \R^d))}
\end{equation}
$$
P(R,Z)\in L^1(I\times\Omega),\; H(R,Z)\in L^1(I;L^1(\Omega))\cap L^1(I; L^1(\partial\Omega; |\vu_B \cdot \vc{n}|
{\rm d}S_x))
$$
is a dissipative turbulent solution to problem \eqref{P1}--\eqref{P7}
iff:
\begin{enumerate}
\item
The integral formulation of the continuity equations
\begin{equation} \label{W1}
\begin{split}
\left[ \intO{ r \varphi } \right]_{t = 0}^{t = \tau} &+
\int_0^\tau \int_{\Gamma^{\rm out}} \varphi r \vu_B \cdot \vc{n} \ \D \ S_x
+
\int_0^\tau \int_{\Gamma^{\rm in}} \varphi r_B \vu_B \cdot \vc{n} \ \D \ S_x\\ &=
\int_0^\tau \intO{ \Big[ r \partial_t \varphi + r \vu \cdot \Grad \varphi \Big] } \dt
\end{split}
\end{equation}
holds
for any $0 \leq \tau \leq T$, and any test function $\varphi \in C^1([0,T] \times \Ov{\Omega})$,
\begin{equation} \label{W2}
r(0, \cdot) = r_0,
\end{equation}
where $r$ stands for $R$ and $Z$.
\item
There exists a tensor measure
\[
 \mathfrak{R} \in L^\infty(0,T; \mathcal{M}^+(\Ov{\Omega}; \R^{d \times d}_{\rm sym})),
\]
such that the integral identity
\begin{equation} \label{W3}
\begin{split}
\left[ \intO{ (R+Z) \vu \cdot \bfphi } \right]_{t=0}^{t = \tau} &=
\int_0^\tau \intO{ \Big[ (R+Z)\vu \cdot \partial_t \bfphi + (R+Z) \vu \otimes \vu : \Grad \bfphi
+ P(R,Z) \Div \bfphi \\&- \mathbb{S}(\Grad \vu) : \Grad \bfphi \Big] }
+ \int_0^\tau \int_{{\Omega}} \Grad \bfphi : \D \ \mathfrak{R}(t) \ \dt
\end{split}
\end{equation}
holds for any $0 \leq \tau \leq T$ and any test function $\bfphi \in C^1([0,T] \times {\Omega}; \R^d)$, $\bfphi|_{\partial \Omega} = 0$,
\begin{equation} \label{W4}
(R+Z) \vu(0, \cdot) = \vm_0:=(R_0+Z_0)\vu_0.
\end{equation}
Here we assume that all quantities appearing in \eqref{W3} are at least integrable in $(0,T) \times \Omega$.

\item
There exists an energy defect measure
\[
\mathfrak{E} \in L^\infty(0,T; \mathcal{M}^+ (\Ov{\Omega}))
\]
such that
\begin{equation} \label{W5}
\begin{split}
&\left[ \intO{\left[ \frac{1}{2} (R+Z) |\vu - \vu_B|^2 + H(R,Z) \right] } \right]_{t = 0}^{ t = \tau} +
\int_0^\tau \intO{ \mathbb{S}(\Grad\vu):\Grad\vu \Big] } \dt\\
&+\int_0^\tau \int_{\Gamma^{\rm out}} H(R,Z)  \vu_B \cdot \vc{n} \ \D S_x \dt +
\int_0^\tau \int_{\Gamma^{\rm in}} H(R_B, Z_B)  \vu_B \cdot \vc{n} \ \D S_x \dt
+ \int_{\Ov{\Omega}}  \D \ \mathfrak{E} (\tau) \\	
\leq
&-
\int_0^\tau \intO{ \left[ (R+Z) \vu \otimes \vu + P(R,Z) \mathbb{I} \right]  :  \Grad \vu_B } \dt + \int_0^\tau { \intO{ {(R+Z)} \vu\cdot \Grad \vu_B   \cdot \vu_B  } }
\dt \\ &+ \int_0^\tau \intO{ \mathbb{S}(\Grad\vu) : \Grad \vu_B } \dt -
\int_0^\tau \int_{\Ov{\Omega}} \Grad \vu_B : \D \ \mathfrak{R}(t) \dt.
\end{split}
\end{equation}
for a.a. $0 \leq \tau \leq T$.

\item Finally, compatibility conditions between
the energy defect $\mathfrak{E}$ and the Reynolds defect $\mathfrak{R}$,
are verified,
\begin{equation} \label{W6}
\underline{d} \mathfrak{E} \leq {\rm Tr}[\mathfrak{R}] \leq \Ov{d} \mathfrak{E},\
\ \mbox{for certain constants}\ 0 < \underline{d} \leq \Ov{d}.
\end{equation}
\end{enumerate}
\end{Definition}

\begin{Remark} \label{CR1}
\begin{enumerate}
\item The { compatibility} condition (\ref{W6}) is absolutely crucial for the weak--strong uniqueness principle stated in Theorem \ref{ET4} below.
\item
In view of \eqref{W6}, one can always consider
\begin{equation} \label{DD1}
\mathfrak{E} \equiv \frac{1}{\Ov{d}} {\rm tr}[\mathfrak{R}];
\end{equation}
whence, strictly speaking, the energy defect $\mathfrak{E}$ can be completely omitted in the definition.
\item
As we shall see in the existence proof below, the dissipative solutions can be constructed in such a way that the constant
$\Ov{d}$ depends solely on the dimension $d$ and the structural constants $\underline{a}$, $\Ov{a}$ appearing in \eqref{Pressure2}.
\item We remark that all conclusions of this paper, after necessary modification of definitions, hold for general
non Newtonian fluids { characterized} by {\it general rheological law} in the spirit of \cite{AbFeNo}. It is also possible to prescribe the Navier boundary conditions (instead of the Dirichlet boundary conditions) on a part of the
$\Gamma^{\rm out}$ boundary.
\end{enumerate}
\end{Remark}

\subsection{Main results}

\subsubsection{Assumptions}
Let
 \begin{equation}\label{calO}
  {\cal O}=\{(R,Z)\,|\, \underline b R<Z<\overline b R\}\; \mbox{with some
  $0<\underline b<\overline b<\infty$.}
 \end{equation}
We suppose that
\begin{enumerate}
 \item {\it Domain:}
 \begin{equation}\label{Omega}
  \Omega\,\mbox{ is a bounded Lipschitz domain.}
 \end{equation}
 \item {\it Boundary data:}
 \begin{equation}\label{Bdata1}
 \vu_B\in C^1_c(\R^d;\R^d), r_B\in C_c(\R^d),\;r_B\ge 0,
 \end{equation}
 where $r_B$ stands for $R_B$, $Z_B$, and
 \begin{equation}\label{Bdata2}
  (R_B,Z_B)\in \overline{\cal O}
 \end{equation}
\item{\it Initial data:} There exists $\gamma>1$ such that
\begin{equation}\label{Idata1}
R_0 \in L^{\gamma} (\Omega),\ R_0 \geq 0,\ \vm_0 \in L^{\frac{2 \gamma}{\gamma + 1}}(\Omega; R^d),
\intO{ \left[ \frac{1}{2} \frac{|\vm_0|^2}{R_0+Z_0} + H(R_0,Z_0) \right] } < \infty,
\end{equation}
and
\begin{equation}\label{Idata2}
  (R_0,Z_0)\in \overline{\cal O}.
 \end{equation}
\item {\it Pressure--density equation of state:}
\begin{equation} \label{Pressure1}
P\in C^1[\overline{\cal O}) \cap C^2({\cal O}),
P(0,0) = 0.
\end{equation}
The Helmholtz function $H$ defined by (\ref{P4}) and $P$ are such that
\begin{equation}\label{Pressure2}
H\;\mbox{is strictly convex  on ${\cal O}$},\;
H - \underline{a} P,\ \Ov{a} P - H\ \mbox{are convex on ${\cal O}$}.
\end{equation}
\end{enumerate}

\begin{Remark}\label{CR2}
\begin{enumerate}
\item An iconic example an equation of state satisfying assumptions \eqref{Pressure1}--
\eqref{Pressure2}
is the isentropic pressure--density relation
\begin{equation}\label{S1}
P(R,Z) = a_1 R^\gamma+a_2 Z^\beta,\ a_1,a_2 > 0, \ \gamma,\beta > 1
\end{equation}
\item One easily checks by using (\ref{P4}), that $H\in C(\overline{\cal O})\cap C^1({\cal O})$
\item We may suppose without loss of generality that also $P$, $H - \underline{a} P$,\ $\Ov{a} P - H$ are strictly convex on ${\cal O}$.
\item Due to \eqref{P4}, $P$ and $H$ are interrelated by the
differential equation
\begin{equation}\label{H+}
 R\partial_RH(R,Z)+ Z\partial_ZH(R,Z)-H(R,Z)=P(R,Z).
\end{equation}
\item
It is easy to check that any $P$ satisfying \eqref{Pressure1}--\eqref{Pressure2} possesses certain coercivity
similar to \eqref{S1}. More specifically,
\begin{equation} \label{Pres2a}
P(R,Z), H(R,Z)\ge a R^\gamma \  \mbox{for all}\ R \geq \overline R,\ (R,Z)\in \overline{\cal O}\, \mbox{for certain}\ a > 0,\ \gamma > 1, \ \overline R>0.
\end{equation}
Indeed as $\Ov{a}P - H$ is a convex function and $H$ is strictly convex, we get
\[
\Ov{a} {\mathfrak{P}}_s''(R) \geq {\mathfrak{ H}}_s''(R) = \frac{\mathfrak{P}_s'(R)}{R},\ R > 0,\, s\in[\underline b,\overline b].
\]
where
$$
\mathfrak{P}_s(R)= P(R,sR),\;\mathfrak{H}_s(R)= H(R,sR).
$$
In particular $\mathfrak{P}'_s(R)>0$ whatever $s\in [\underline b,\overline b]$ is, and by the uniform continuity $\inf_{s\in  [
\underline b,\overline b]}\mathfrak{P}'_s(1)=\underline c>0$.
Moreover, since $\mathfrak{P}_s(0)=0$ we also have $P\ge 0$.

This yields
\[
{ \Big(\log(\mathfrak{P}_s'(R))\Big)' \geq \Big(\log \left(R^{\frac{1}{\Ov{a}}} \right)\Big)'} \ \Rightarrow \
P(R,Z) \geq\underline c R^{1+\frac{1}{\Ov{a}}} \ \mbox{for all}\ R\ge \overline R, \
(R,Z)\in {\cal O},
\]
and consequently, since $\mathfrak{H}''_s(R)=\mathfrak{P}'_s(R)/R$,
$$
H(R,Z) \geq\underline c R^{1+\frac{1}{\Ov{a}}} \ \mbox{for all}\ R\ge \overline R, \
(R,Z)\in {\cal O},
$$
whence \eqref{Pres2a} holds for $\gamma = 1 + \frac{1}{\Ov{a}}$.
\item There exist non negative real numbers  $\mathfrak{a}_1$, $\mathfrak{a}_2$, $\mathfrak{a_3}$ such that
$$
0\le P(R,Z)\le H(R,Z)+\mathfrak{a}_1 R+\mathfrak{a}_2 Z+\mathfrak{a_3}
$$

\end{enumerate}
\end{Remark}

We are now able to formulate {the} main results of this paper.

\subsubsection{Existence}

\begin{Theorem} [{\bf Global existence of dissipative turbulent solutions}] \label{ET1}
Let assumptions \eqref{Omega}--\eqref{Pressure2} be satisfied.
Then the problem \eqref{P1}--\eqref{P7} admits at least one dissipative turbulent solution $[R,Z,\vu]$ in $(0,T) \times \Omega$ in the sense specified in Definition \ref{DW1}.
\end{Theorem}

\subsubsection{Compatibility with classical solutions}
\begin{Theorem} [{\bf Compatibility of regular turbulent solutions with classical solutions}] \label{ET2}
Let assumptions \eqref{Omega} and \eqref{Pressure1}
be satisfied.
Suppose that $[R,Z,\vu]$ is a dissipative turbulent solution to problem
(\ref{P1})--(\ref{P7}) in the sense of Definition \ref{DW1} belonging to the class
$$
\vu\in C^1(\overline I\times\overline\Omega;\R^d),\;R,Z \in C^1(\overline I\times\overline\Omega),\;
\inf_{I\times\Omega}Z>0.
$$
Then $\mathfrak{E}=\mathfrak{R}=0$ and equations (\ref{P1})--(\ref{P7})
are satisfied in the classical sense.
\end{Theorem}

\subsubsection{Relative energy inequality}

We introduce the relative energy functional
\begin{equation}\label{calE}
\mathcal{E} \left(R, Z, \vu \ \Big|\ \mathfrak{r}, \mathfrak{z},\mathfrak u \right) = \frac{1}{2} (R+Z)|\vu - \mathfrak{u}|^2 + H(R,Z) -
\partial_RH(\mathfrak{r},\mathfrak{z}) (R - \mathfrak{r})
-\partial_ZH(\mathfrak{r},\mathfrak{z})(Z-\mathfrak{z})- H(\mathfrak{r},\mathfrak{z}),
\end{equation}
where $R,Z,\vu$ is a dissipative turbulent solution of problem
\eqref{P1}--\eqref{P7}, while $\mathfrak{r},\mathfrak{z},\mathfrak{u}$
are test functions in class:
$$
\mathfrak{u}\in C^1(\overline I\times\overline\Omega;\R^d),\;
\Div\mathbb{S}(\Grad\mathfrak{u})\in C(\overline I\times\overline\Omega;\R ^d),\; \vu|_{\partial\Omega}=\vu_B|_{\partial\Omega},
$$
\begin{equation}\label{RR6}
\inf_{I\times\Omega}\mathfrak{z}>0,\; (\mathfrak{r},\mathfrak{z})\in C^1(\overline I\times\overline\Omega;\overline{\cal O}).
\end{equation}

The following Theorem describes the evolution of ${\cal E}$:
\begin{Theorem}[{\bf Relative energy inequality}] \label{ET3}
Suppose that $\Omega$ is a bounded
Lipschitz domain and that $P$ satisfies hypotheses (\ref{Pressure1}).
Let $[R, Z,\vu]$ be a dissipative turbulent solution in the sense of Definition \ref{DW1}.  Then there holds:
\begin{equation}\label{RR5}
\begin{split}
&\left[ \intO{\mathcal{E}\left(R,Z, \vu \ \Big|\ \mathfrak{r}, \mathfrak{z}, \mathfrak{u} \right) } \right]_{t = 0}^{ t = \tau} +
\int_0^\tau \intO{\mathbb S(\Grad\vu):\Grad(\vu- \mathfrak{u}) } \dt
\\
&+\int_0^\tau \int_{\Gamma^{\rm out}} \left[ H(R,Z) - \partial_RH(\mathfrak{r}, \mathfrak{z}) (R - \mathfrak{r}) -
\partial_ZH(\mathfrak{r}, \mathfrak{z}) (Z - \mathfrak{z})-
H(\mathfrak{r})  \right]  \vu_B \cdot \vc{n} \ \D S_x \dt
\\
&+\int_0^\tau \int_{\Gamma^{\rm in}} \left[ H(R_B,Z_B) - \partial_RH(\mathfrak{r}, \mathfrak{z}) (R_B - \mathfrak{r}) -
\partial_ZH(\mathfrak{r}, \mathfrak{z}) (Z_B - \mathfrak{z})-
H(\mathfrak{r})  \right]  \vu_B \cdot \vc{n} \ \D S_x \dt
\\	
&+ \int_{\Ov{\Omega}} 1 \D \ \mathfrak{E} (\tau)\le - \int_0^\tau \intO{ (R+Z) (\mathfrak{u} - \vu) \cdot\Grad \mathfrak{u}\cdot (\mathfrak{u} - \vu)  } \dt
\\
&-
\int_0^\tau \intO{ \left[ P(R,Z) - \partial_RP(\mathfrak{r}, \mathfrak{z}) (R - \mathfrak{r}) -
\partial_ZP(\mathfrak{r}, \mathfrak{z}) (Z - \mathfrak{z})-
P(\mathfrak{r},\mathfrak{z})  \right] \Div \mathfrak{u} } \dt
\\
&+ \int_0^\tau \intO{ \Big(\frac{R+Z}{\mathfrak{r}+\mathfrak{z}} -1\Big)(\mathfrak{u} - \vu) \cdot \Big[ \partial_t ((\mathfrak{r}+\mathfrak{z}) \mathfrak{u})  +  \Div ((\mathfrak{r}+\mathfrak{z}) \mathfrak{u} \otimes \mathfrak{u})
   \Big] } \dt
  \\
&+ \int_0^\tau \intO{ (\mathfrak{u} - \vu) \cdot \Big[ \partial_t ((\mathfrak{r}+\mathfrak{z}) \mathfrak{u})  +  \Div ((\mathfrak{r}+\mathfrak{z}) \mathfrak{u} \otimes \mathfrak{u})
  + \Grad P(\mathfrak{r},\mathfrak{z}) \Big] } \dt
  \\
  &+ \int_0^\tau \intO{
(R-\mathfrak{r})(\mathfrak{u}-\vu)\cdot\Big(\nabla\mathfrak{r}\partial^2_RH(\mathfrak{r},\mathfrak{z})
+\nabla\mathfrak{z}\partial_R\partial_Z H(\mathfrak{r},\mathfrak{z})\Big)}{\rm d}t
\\
&+\int_0^\tau \intO{
(Z-\mathfrak{z})(\mathfrak{u}-\vu)\cdot\Big(\nabla\mathfrak{z}\partial^2_ZH(\mathfrak{r},\mathfrak{z})
+\nabla\mathfrak{r}\partial_R\partial_Z H(\mathfrak{r},\mathfrak{z})\Big)
  } \dt
  \\
  &+ \int_0^\tau \intO{\Big[\frac{Z+R}{\mathfrak{r}+\mathfrak{z}}
  (\vu-\mathfrak{u})\cdot\mathfrak{u}+  \partial_RP(\mathfrak{r},\mathfrak{z})-R\partial^2_RH(\mathfrak{r},\mathfrak{z})-Z\partial_R\partial_ZH(\mathfrak{r},\mathfrak{z
  })\Big]\Big[ \partial_t \mathfrak{r}   +  \Div (\mathfrak{r} \mathfrak{u})
  \Big] } \dt
  \\
  &+\int_0^\tau \intO{\Big[\frac{Z+R}{\mathfrak{r}+\mathfrak{z}}
  (\vu-\mathfrak{u})\cdot\mathfrak{u}+\partial_ZP(\mathfrak{r},\mathfrak{z})-  { Z\partial^2_Z}H(\mathfrak{r},\mathfrak{z})-{ R\partial_R}\partial_ZH(\mathfrak{r},\mathfrak{z
  })\Big] \Big[ \partial_t \mathfrak{z}   +  \Div (\mathfrak{z} \mathfrak{u})
  \Big] } \dt
  \\
&- \int_0^\tau \int_{\Ov{\Omega}} \Grad \mathfrak{u} : \D \ \mathfrak{R}(t) \dt
\end{split}
\end{equation}
with any $(\mathfrak{r},\mathfrak{z},\mathfrak{u})$ in class (\ref{RR6}).
\end{Theorem}

\subsubsection{Weak-strong uniqueness}
\begin{Theorem} [{\bf Weak--strong uniqueness}] \label{ET4}
Let assumptions \eqref{Omega}--\eqref{Pressure2} be satisfied.
Let $[R,Z, \vu]$ be a dissipative solution
in the sense of Definition \ref{DW1}, and let $[\mathfrak{r},\mathfrak{z}, \mathfrak{u}]$ be a strong solution of the same problem belonging to the
class (\ref{RR6})
Then
\[
R = \mathfrak{r}, \ Z=\mathfrak{z},\ \vu = \mathfrak{u} \ \mbox{in}\ (0,T) \times \Omega,\ \mathfrak{E} = \mathfrak{R} = 0.
\]
\end{Theorem}

The following remark to Theorems \ref{ET1}--\ref{ET4} is in order:
\begin{Remark}\label{CR3}
\begin{enumerate}
 \item The value of $\gamma$ in Theorem \ref{ET1} (cf. Definition \ref{DW1}) is the minimum of
$\gamma$ from assumption (\ref{Idata1}) and $\gamma$ calculated
in Remark \ref{CR2}, cf. formula \eqref{Pres2a}.
\item We notice that compatibility theorem (Theorem \ref{ET2})
as well as relative energy inequality (Theorem \ref{ET3})
do not require practically any structural assumptions on the pressure.
\item
It is to be noticed that the isothermal pressure $P(R,Z)=a_1R+a_2Z$, $a_i>0$ does not satisfy the hypothhesis (\ref{Pressure2}). These conditions are however necessary
for the Reynolds stress $\mathfrak{R}$ to be a positively semi-definite tensor.
Thus, from the point of view of physics, conditions \eqref{Pressure2} may seem too restrictive. A brief inspection of the proofs reveals that all principal results remain valid for any equation of state of the  form
\[
P(R,Z) + a_1R+a_2Z, \ a_i \geq 0
\]
as long as $P$ satisfies \eqref{Pressure2}. To see it, one has to take advantage of the linearity
of the ``perturbation'' $a_1R+a_2Z$ in the limiting process in the proofs.
\end{enumerate}

\end{Remark}

\subsubsection{A remark on local existence on strong solutions}

Theorem \ref{ET4} operates with strong solutions to problem \eqref{P1}--\eqref{P7}. A question of existence of such solutions at least locally in time is therefore quite natural. Such results are however in a short supply even for a slightly more simple monofluid case. To the best of our knowledge, all of them require quite particular geometrical
conditions on the inflow boundary. One of the most representative sample of such results is Theorem 2.5
in Valli and Zajaczkowski \cite{VAZA}. Its
reformulation to the bi-fluid system \eqref{P1}--\eqref{P7} leads to the following statement (compare with \cite{JiNo}, where the author treat the case of zero inflow-outflow):
\\

{\it
Let $\Omega\in C^3$ be a bounded domain, $0< \underline{\mathfrak{r}}<\overline{\mathfrak{r}}<\infty$,
$0< \underline{\mathfrak{z}}<\overline{\mathfrak{z}}<\infty$,
be constants. Assume that
$$
P\in C^2((0,\infty)^2).
$$
Suppose that
$$
\vu_B\in W^{3,2}(\Omega),\;\mathfrak{r}_B\in W^{2,2}(\Omega),
$$
$$
\vu_B\cdot\vc n|_{\Gamma^{\rm in}}\ge \underline u>0,
$$
$$
\vu_0-\vu_B\in W^{1,2}_0(\Omega), \; \mathfrak{r}_0, \mathfrak{z}_0\in W^{2,2}(\Omega),
$$
$$
 \underline{\mathfrak{ r}}\le\mathfrak{r}_0\le\overline{\mathfrak{r}}_0,\;
\; \underline {\mathfrak{z}}\le\mathfrak{z}_0\le\overline{\mathfrak{z}},
$$
$$
\Div(\vr_0\vu_0)|_{\Gamma^{\rm in}}=0,
$$
$$
\frac 1{\mathfrak{r}_0+\mathfrak{z}_0}\Big(-\nabla P(\mathfrak{r}_0,\mathfrak{z}_0)+\mu \Delta\vu_0
+(\mu +\lambda)\nabla{\rm div}\vu_0 -(\mathfrak{r}_0+\mathfrak{z}_0)\vu_0\nabla\vu_0\Big)
\in W^{1,2}_0(\Omega).
$$
\begin{enumerate}
\item
Then there exists an interval $I_*=[0,T_*)$ and numbers $\underline r$, $\overline r$, $\underline z$, $\overline z$, $0< \underline r<\underline{\mathfrak{r}}<\overline{\mathfrak{r}}<\overline r<\infty$,
$0< \underline z<\underline{\mathfrak{z}}<\overline{\mathfrak{z}}<\overline z<\infty$
such that the problem (\ref{P1}--\ref{P7}) admits in the class
\begin{equation}\label{st1}
(\mathfrak{r},\mathfrak{z})\in C(I_*;W^{2,2}(\Omega)),\,\partial_t(\mathfrak{r},\mathfrak{z})\in C(I_*;W^{1,2}(\Omega)),
\end{equation}
$$
\vu\in L^2(I_*;W^{3,2}(\Omega;\R^3)),\,\partial_t\vu\in L^2(I_*;W^{2,2}(\Omega,\R^3)),\;\partial^2_t\vu\in L^2(I_*;L^{2}(\Omega,\R^3)),
$$
\begin{equation}\label{stm2}
\underline r\le \mathfrak{r}\le\overline r,\, \underline z\le \mathfrak{z}\le\overline z,
\end{equation}
$$
 \vu(0)=\vu_0,\;\vu|_{(0,T)\times\partial\Omega}=\vu_B|_{\partial\Omega},\;\mathfrak{r}|_{\Gamma^{\rm in}}=
 \mathfrak{r}_B|_{\Gamma^{\rm in}}
 $$
a unique strong solution
$(\mathfrak{r}, \mathfrak{z},\vu)$.

\item
If moreover
$$
 \underline b \mathfrak{r}_B|_{\Gamma^{\rm in}}\le\mathfrak{z}_B|_{\Gamma^{\rm in}}\le\overline b \mathfrak{r}_B|_{\Gamma^{\rm in}},\;
 \underline b \mathfrak{r}_0\le\mathfrak{z}_0\le\overline b \mathfrak{r}_0,
$$
with some
$0<\underline b<\overline b<\infty$, then
\begin{equation}\label{stm3}
\underline b \mathfrak{r}\le \mathfrak{z}\le\overline b \mathfrak{r}.
\end{equation}
\end{enumerate}
}

Condition $\vu_B\cdot\vc n|_{\Gamma^{\rm in}}\ge\underline u>0$
is very restrictive. In practical situations, it can be satisfied only provided $\Gamma^{\rm in}$ is a union of nonintersecting compact manifolds.

Local existence of strong solutions with non-zero inflow/outflow in general case is, even for the mono-fluid models (with one continuity equation), to our best knowledge, an open problem. The essence of the difficulties dwells in the conditions allowing sufficiently smooth extensions of velocity field outside $\overline\Omega$ and in the the ``management'' of flow corresponding to the extended velocity field. These difficulties are the same in the mono-fluid case
(they are independent on the number of continuity/transport equations in the system that need to be treated).

\section{Existence (Proof of Theorem \ref{ET1})}
\label{E}

Our first goal is to show that the dissipative solutions exist globally in time for any finite energy initial data.

The proof is based on a multilevel approximation scheme that shares certain features with the approximation of the
compressible Navier--Stokes equations in \cite{ChJiNo}, \cite{AbFeNo}.

\subsection{First level approximation}

First, we introduce a sequence of finite--dimensional spaces
$X_n \subset L^2(\Omega; \R^d)$,
\[
X_n = {\rm span} \left\{ \vc{w}_i\ \Big|\ \vc{w}_i \in \DC(\Omega; \R^d),\ i = 1,\dots, n \right\}
\]
Without loss of generality, we may assume that $\vc{w}_i$ are orthonormal with respect to the standard scalar product in $L^2$. We denote by $\Pi_n$ the orthogonal projection
of $L^2(\Omega)$ onto $X_n$.

Concerning initial data, we  may suppose without loss of generality that initial and boundary data are smooth and strictly positive, i.e. in addition to (\ref{Bdata2}), (\ref{Idata2}),
\begin{equation}\label{Bbis}
r_B\in C^1(\R^d),\,0<\underline r\le r_B\le \overline r<\infty,\; \vu_B\in C_c^1(\R^d;\R^d),
\end{equation}
\begin{equation}\label{Ibis}
r_0\in C^1(\R^d),\,0<\underline r\le r_0\le \overline r<\infty,\;\vu_0\in C_c^1(\R^d;\R^d).
\end{equation}
In \eqref{Bbis}--\eqref{Ibis} $r_0$, $r_B$ stands for $R_0$, $R_B$ and $Z_0$, $Z_B$,
respectively.

Following \cite{ChJiNo},
 we use a parabolic approximation of the equations of continuity,
\begin{equation} \label{E1}
\partial_t r + \Div (r \vu ) = \ep \Del r \ \mbox{in}\ (0,T) \times \Omega,\ \ep > 0,
\end{equation}
supplemented with the boundary conditions
\begin{equation} \label{E2}
\ep \Grad r \cdot \vc{n} + (r_B - r) [\vu_B \cdot \vc{n}]^- = 0 \ \mbox{in}\ [0,T] \times \partial \Omega,
\end{equation}
and the initial condition
\begin{equation} \label{E3}
r(0, \cdot) = r_0.
\end{equation}
Here $r$ stand for $R$ and $Z$ and $\vu = \vv + \vu_B$, with $\vv \in C([0,T]; X_n)$, in particular, $\vu|_{\partial \Omega} = \vu_B$, { and symbol $[a]^-:=\min\{a,0\}$}.
Note that for given $\vu$, $r_B$, $\vu_B$, this is a linear problem for the unknown $r$.

As $\Omega$ is merely Lipschitz, the usual parabolic estimates fail at the level of the spatial derivatives (we cannot use at this stage
the maximal regularity theory as it was done in \cite{ChJiNo})
and we
are forced to use the weak formulation:
\begin{equation} \label{E4}
\begin{split}
\left[ \intO{ r \varphi } \right]_{t = 0}^{t = \tau}&=
\int_0^\tau \intO{ \left[ r \partial_t \varphi + r \vu \cdot \Grad \varphi - \ep \Grad r \cdot \Grad \varphi \right] }
\dt \\
&- \int_0^\tau \int_{\partial \Omega} \varphi r \vu_B \cdot \vc{n} \ {\rm d} S_x \dt +
\int_0^\tau \int_{\partial \Omega} \varphi (r - r_B) [\vu_B \cdot \vc{n}]^{-}  \ {\rm d}S_x \dt,\
r(0, \cdot) = r_0,
\end{split}
\end{equation}
for any test function
\[
\varphi \in L^2(0,T; W^{1,2}(\Omega)),\ \partial_t \varphi \in L^1(0,T; L^2(\Omega))
\]
as in \cite{AbFeNo}.

Following \cite{KwNo}, \cite{NoSCM}, we use the Galerkin apparoximation to approximate the momentum equation: We look for approximate velocity field in the form
\[
\vu = \vv + \vu_B, \ \vv \in C([0,T]; X_n).
\]
Accordingly, the
approximate momentum balance reads
{
\begin{equation} \label{E5}
  \intO{ (R+Z) \vu \cdot \bfphi } \Big|_{t=0}^{t = \tau} =
\int_0^\tau \intO{ \Big[ (R+Z) \vu \cdot \partial_t \bfphi + (R+Z) \vu \otimes \vu : \Grad \bfphi
+ P(R,Z) \Div \bfphi
\end{equation}
$$
- \mathbb{S}(\nabla\vu): \Grad \bfphi -\ep\Grad (R+Z) \cdot \Grad \vu \cdot \bfphi \Big] }\dt
$$
}
%
for any $\bfphi \in C^1([0,T]; X_n)$, with the initial condition
\begin{equation} \label{E6}
(R+Z) \vu(0, \cdot) = (R_0+Z_0) \vu_0, \ \vu_0 = \vv_0 + \vu_B, \ \vv_0=\Pi_n(\vu_0-\vu_B).
\end{equation}
For fixed parameters $n$, $\ep > 0$, the first level approximation is a solution $[R, Z,  \vu]$
\footnote{Here in the sequel, we skip the indexes  $\ep$, $n$ and write e.g. $R$ instead of
$R_{\ep,n}$, etc. and will use eventually only one of them in the situations when it will be useful to underline the corresponding limit passage.}
of the parabolic problem \eqref{E1}--\eqref{E3}, and the Galerkin approximation \eqref{E5}, \eqref{E6}.

\subsection{Parabolic problem with Robin boundary conditions}

In setting (\ref{E4}) on Lipschitz domains (and even in a more general setting as far as the regularity of the transporting velocity $\vu$ is concerned) problem
\eqref{E1}--\eqref{E3} has been investigated in Crippa, Donadello, Spinolo \cite{CrDoSp}. The following lemma resumes
Lemmas \cite[Lemma 3.2 and Lemma 3.4]{CrDoSp} and \cite[Lemma 3.3,
Corollary 3.4 and estimate (3.7)]{AbFeNo}:

\begin{Lemma} \label{EL1}
Let $\Omega \subset \R^3$ be a bounded Lipschitz domain and $\vu = \vv + \vu_B$, $\vv \in C(\overline I; X_n)$. Suppose that
$(r_B,\vu_B)$ belongs to the class (\ref{Bbis}) while $r_0$ belongs to the class (\ref{Ibis}). Then we have:
\begin{enumerate}
\item
The initial--boundary value problem (\ref{E1}--\ref{E3})
admits a weak solution $r$ {specified in (\ref{E4})}, unique in the class
$$
r \in L^2(I; W^{1,2}(\Omega)) \cap C(\overline I; L^2(\Omega)).
$$
The norm in the aforementioned spaces is bounded only in terms of the data $r_B$, $r_0$, $\vu_B$, and
$\|\vv, {\rm div} \vv\|_{L^\infty(I;L^\infty(\Omega))}$.
\item Moreover, $\partial_t r\in L^2(I\times\Omega)$ and $\sqrt\ep\nabla r\in L^\infty(I;L^2(\Omega))$ are bounded in terms of the data $r_B$, $r_0$, $\vu_B$ and
$\|\vv, {\rm div} \vv\|_{L^\infty(I;L^\infty(\Omega))}$ and $\nabla^2 r\in L^2(I; L^2_{\rm loc}(\Omega)$ is bounded in the same way on any compact set $K$ of $\Omega$
with the constant dependent in addition on $K$.
\item Strong maximum principle: The solution satisfies,
\begin{equation}\label{max}
\begin{array}{c}
{
\forall \tau\in\overline I,\;\| r(\tau) \|_{L^\infty( \Omega)}
}
 \leq M
 \exp \left( T \| \Div \vu \|_{L^\infty((0,\tau) \times \Omega)} \right),
 \\ \\
\mbox{for a.a. $\tau\in I$},\; r(\tau,x)\leq M\exp \left( T \| \Div \vu \|_{L^\infty((0,\tau) \times \Omega)} \right)\;\mbox{for a.a. $x\in \partial\Omega$},
\end{array}
\end{equation}
where
$$
M=\max \left\{\max_\Omega r_0,\max_{\Gamma^{\rm in}} r_B,
\| \vu_B \|_{L^\infty((0,T) \times \Omega)} \right\}.
$$
\item Strong minimum principle: The solution satisfies,
\begin{equation}\label{min}
\begin{array}{c}
 \forall \tau\in\overline I,\;{\rm ess} \inf_{x\in \Omega} r(\tau,x) \geq
m
\exp \left( -T \| \Div \vu \|_{L^\infty((0,T) \times \Omega)} \right),
\\ \\
\mbox{for a.a. $\tau\in I$},\; r(\tau,x)\ge m\exp \left(- T \| \Div \vu \|_{L^\infty((0,\tau) \times \Omega)} \right)\;\mbox{for a.a. $x\in \partial\Omega$},
\end{array}
\end{equation}
where
$$
m=\min \left\{ \min_\Omega r_0 , \min_{\Gamma^{\rm in}} r_B \right\}.
$$
\end{enumerate}
\end{Lemma}

\subsection{Existence of first level approximation { ($\ep$, $n$ fixed)}}

The \emph{existence} of the approximate solutions at the level of
the parabolic problem (\ref{E1}--\ref{E3}) coupled with the Galerkin approximation (\ref{E5}--\ref{E6})
can be proved in the same way as in \cite[Section 4]{ChJiNo} (mono-fluid case with non zero inflow-outflow) combined with \cite[Section 3]{NoSCM},
eventually with \cite[Section 4]{AN-MP} (multi-fluid with zero boundary conditions).  Specifically,
for $\vu = \vu_B + \vv$, $\vv \in C([0,T]; X_n)$, we identify the unique solutions $r =
r [\vu]$ of (\ref{E1}--\ref{E3}), where $r$ stands for $R, Z$ and plug them as $R$, $Z$ in \eqref{E5}.
The unique solution $\vu = \vu[R,Z]$ of \eqref{E5} defines a mapping
\[
\mathcal{T}: \vv \in C([0,T]; X_n) \mapsto   \mathcal{T}[\vc{v}] = (\vu[R,Z] - \vu_B)
\in C([0,T]; X_n).
\]
The first level approximate solutions $r = r_{n,\ep}$, $\vu = \vu_{n,\ep}$ -- here, $r$ stands for $R$, $Z$--are obtained via a fixed point
of the mapping $\mathcal{T}$. This  procedure is detailed in \cite{ChJiNo} and in { \cite{KwNo}} for the mono-fluid case
with the non zero inflow-outflow and in \cite{NoSCM} for the multi-fluid case with the no-slip boundary conditions. Combinig
\cite[Section 4]{KwNo} with \cite[Section 4]{NoSCM},
we easily deduce the following result.\footnote{ The energy inequality (\ref{E7}) in \cite[Lemma 4.2]{KwNo} and in \cite[Section 4]{NoSCM}
is derived under assumption $\Omega\in C^2$. This assumption is needed due to the treatment of the parabolic problem (\ref{E1}--\ref{E3}) via
the classical maximal regularity methods. With Lemma \ref{EL1} at hand, the same proof can be carried out
without modifications also in Lipschitz domains.}

\begin{Proposition}[{\rm { First level approximate solutions ($\ep$, $n$ fixed)}}] \label{EP1}

Let $\Omega \subset \R^3$ be a bounded Lipschitz domain.
Let the data $(R_B,Z_B,\vu_B)$, $(R_0,Z_0,\vu_0)$  belong to the class { (\ref{Bdata2}), (\ref{Bbis})},
(\ref{Idata2}), (\ref{Ibis}). Suppose that assumptions (\ref{Pressure1}--\ref{Pressure2}) are satisfied.

Then for each fixed $n > 0$, $\ep > 0$,  there exists a solution
$(R_\ep,Z_\ep,\vu_\ep=\vv_\ep+\vu_B)$
in the class
$$
R,Z\in L^2(I;W^{1,2}(\Omega))\cap L^\infty( I\times\Omega),\
\partial_t(R,Z)\in L^2(I\times\Omega),\ R,Z\in L^\infty(I\times\partial\Omega),
$$
$$
\forall t\in\overline I,\
R(t),Z(t)> 0\ \mbox{a.e. in $\Omega$},\;
\mbox{for a.e. $t\in I$,}\
R(t),Z(t)> 0\ \mbox{a.e. in $\partial\Omega$},
$$
$$
\vv=C(\overline I; X_n), \;\partial_t\vv\in L^2(I;X_n).
$$
of the approximate problem (\ref{E4}) and \eqref{E5}, \eqref{E6}.
Moreover, the following holds:
\begin{enumerate}
\item Lower and upper bounds of "densities":
{
\begin{equation}\label{eq5.1}
\begin{array}{c}
\forall t\in\overline I,\;
0<\underline c(n)\le
R_\ep(t,x),Z_\ep(t,x)\le \overline c(n),\; \underline b R_\ep(t,x) \leq Z_\ep(t,x) \leq \overline b R_\ep(t,x),
\mbox{a.e. in $\Omega$},
\\
\mbox{for a.a. $t\in I$},\;
0<\underline c(n)\le R_\ep(t,x), Z_\ep(t,x) \le\overline c(n),\;\; \underline b R_\ep(t,x) \leq Z_\ep(t,x) \leq \overline b R_\ep(t,x),
\mbox{a.e. in $\partial\Omega$},
\end{array}
\end{equation}
}
\item The approximate energy inequality
\begin{equation} \label{E7}
\begin{split}
&\left[ \intO{\left[ \frac{1}{2} (R_\ep+Z_\ep)|\vv_\ep|^2 + { H}(R_\ep,Z_\ep)  \right] } \right]_{t = 0}^{ t = \tau} +
\int_0^\tau \intO{\mathbb{S}(\Grad\vu_\ep):\Grad\vu_\ep } \dt
\\
&+\int_0^\tau \int_{\Gamma^{\rm out}} { H}(R_\ep,Z_\ep)  \vu_B \cdot \vc{n} \ {\rm d} S_x \dt
- \int_0^\tau \int_{\Gamma^{\rm in}} E_{{ H}}(R_B,Z_B|R_\ep,Z_\ep) \vu_B \cdot \vc{n}
\ {\rm d}S_x \ \dt \\
&+
\ep \int_0^\tau \intO{ { \nabla_{R,Z}^2 {H}(R_\ep,Z_\ep)} [\Grad R_\ep,\Grad Z_\ep] } \dt \\
&\leq
-
\int_0^\tau \intO{ \left[ (R_\ep+z_\ep) \vu_\ep \otimes \vu_\ep + P(R_\ep,Z_\ep) \mathbb{I} \right]  :  \Grad \vu_B } \dt \\
&+ \int_0^\tau { \intO{ (R_\ep+Z_\ep) \vu_\ep  \cdot\Grad \vu_B \cdot  \vu_B  } }
\dt
+ \int_0^\tau \intO{ \mathbb{S}(\Grad\vu_\ep) : \Grad \vu_B } \dt
\\
&- \int_0^\tau \int_{\Gamma^{\rm in}} { H}(R_B,Z_B)  \vu_B \cdot \vc{n} \ {\rm d} S_x \dt
\end{split}
\end{equation}
holds for any $0 \leq \tau \leq T$, where
$$
\nabla_{R,Z}^2 {H}_\delta(R,Z) [\Grad R,\Grad Z]=
\partial^2_RH(R,Z)|\nabla R|^2+2\partial_{R}\partial_ZH\nabla R\cdot\nabla Z
+\partial^2_ZH(R,Z)|\nabla Z|^2.
$$
\end{enumerate}
\end{Proposition}

This is level I of approximations (with two parameters $n$, $\ep$). We shall pass first to the limit $\ep\to 0$ in order to obtain level II of approximations (with one parameters $\ep$). Then we obtain the dissipative turbulent solutions of problem (\ref{P1}--\ref{P7}) by letting $n\to \infty$.

\subsection{The second level approximation {(limit $\ep\to 0$)}}

Our next goal is to send $\ep \to 0$ in the viscous approximation \eqref{E4}, \eqref{E5}, \eqref{E7} for $n$ fixed. In what follows { $\ep \to 0$} mean
limit over a conveniently chosen subsequence (relabeling is not indicated).

\subsubsection{Limit in the { approximate } continuity equations}

 Seeing that $X_n$ is a finite dimensional normed space  and that $H$ is strictly convex, we deduce from (\ref{eq5.1})
 and (\ref{E7}), in particular,
\begin{equation} \label{E13}
\| \vue \|_{L^\infty(0,T; W^{1,\infty}(\Omega))} \leq c,
\end{equation}
\begin{equation} \label{E14}
\mbox{for all $t\in\overline I$},\;
0 < \underline{c} \leq r_\ep(t,x) \leq \Ov{c}\,\mbox{and}
\, (R_\ep(t,x),Z_\ep(t,x))\in \overline{\cal O}\ \mbox{for a.a. $x\in \Omega$},
\end{equation}
$$
\mbox{for a.a. $t\in\overline I$},\;
0 < \underline{c} \leq r_\ep(t,x) \leq \Ov{c}\,\mbox{and}
\, (R_\ep(t,x),Z_\ep(t,x))\in \overline{\cal O} \ \mbox{for a.a. $x\in \partial\Omega$},
$$
\begin{equation} \label{E15}
\ep \left\| \Grad r_\ep\right\|^2_{L^2(I\times\Omega; \R^d)} \leq c.
\end{equation}
In the above and in what follows, $r_\ep$ stands for $R_\ep$ and $Z_\ep$

In view of the uniform bounds established above, we may assume
\begin{equation}\label{E15+}
r_\ep \to r \ \mbox{weakly-(*) in}\ L^\infty((0,T) \times \Omega)
\ \mbox{and weakly in}\ C_{\rm weak}([0,T]; L^r(\Omega)) \ \mbox{for any}\ 1 < r < \infty,
\end{equation}
passing to a suitable subsequence as the case may be. Note that the second convergence follows from the bound on the
time derivative $\partial_t r_\ep$ obtained from equation \eqref{E4}, via an Arzela-Ascoli type compactness argument. We also have
\begin{equation}\label{E15++}
r_\ep \to r \ \mbox{weakly-(*) in}\ L^\infty((0,T) \times \partial \Omega; \D S_x).
\end{equation}
In addition, the limit density admits the same upper and lower bounds as in \eqref{E14}.

Similarly,
\begin{equation} \label{E15b}
\vue \to \vu \ \mbox{weakly-(*) in}\ L^\infty(0,T; W^{1,\infty}(\Omega; \R^d)),
\end{equation}
and
\[
(R_\ep+Z_\ep) \vue \to \vc{m} \ \mbox{weakly-(*) in} \ L^\infty((0,T) \times \Omega; \R^d)).
\]
Moreover, an abstract version of Arzela--Ascoli theorem yields
\begin{equation} \label{E15c}
\vm = (R+Z) \vu \ \mbox{a.a. in}\ (0,T) \times \Omega.
\end{equation}

This is enough to pass to the limi $\ep\to 0$ in the parabolic problem
(\ref{E4}) and to obtain

\begin{equation} \label{E16}
\begin{split}
\left[ \intO{ r \varphi } \right]_{t = 0}^{t = \tau}&=
\int_0^\tau \intO{ \Big[ r \partial_t \varphi + r \vu \cdot \Grad \varphi  \Big] }
\dt \\
&- \int_0^\tau \int_{\Gamma^{\rm out}} \varphi r \vu_B \cdot \vc{n} \ {\rm d} S_x \dt -
\int_0^\tau \int_{\Gamma^{\rm in}} \varphi r_B \vu_B \cdot \vc{n}  \ {\rm d}S_x \dt,\
r(0, \cdot) = r_0
\end{split}
\end{equation}
for any $\varphi \in C^1([0,T] \times \Ov{\Omega})$,
which is a weak formulation of the equation of continuity \eqref{P1}, with the boundary conditions \eqref{P5}, and the initial condition
\eqref{P7}.
\subsubsection{Limit in the approximate momentum equation}

Clearly,
$$
\|\Grad\vue\|_{L^2(I\times\Omega)}\le c
$$
and
$$
\Grad\vue\to\Grad\vu\ \mbox{weakly in {  $L^2(I\times\Omega)$}}.
$$

Next, we deduce from \eqref{E5} on one hand and from Item 2. of Lemma \ref{EL1} on the other hand that
\[
\partial_t
\Pi_n[ (R_\ep+Z_\ep) \vue ] \ \mbox{bounded in}\ L^2(0,T; X_n), \;
\partial_t r_\ep\ \mbox{bounded in}\ L^2(I\times\Omega),
\]
where $\Pi_n: L^2 \to X_n$ is the associated orthogonal projection;
consequently
$$
\|\partial_t\vu_\ep\|_{L^2(I;X_n)}\le c
$$
and, due to Arzela-Ascoli (or Lions-Aubin) compactness argument, we may assume that
$$
\vu_\ep\to\vu\;\mbox{ in}\; C(\overline{ I\times\Omega});
$$
\begin{equation}\label{dod0}
r_\ep\vu_\ep \to r \vu,\
r_\ep \vue \otimes \vue \to r \vu \otimes \vu \
\ C_{\rm weak}(\overline I,L^q(\Omega)),\ 1\le q<\infty.
\end{equation}

By virtue of (\ref{E7}) and the last item in Remark \ref{CR2},
$$
\sup_{\tau\in\overline I}
\|H(R_\ep, Z_\ep)\|_{L^\infty(\Omega)},\,\sup_{\tau\in\overline I}
\|P(R_\ep, Z_\ep)\|_{L^\infty(\Omega)},\;\mbox{is bounded uniformly with
$\ep$, $n$}.
$$

Thus, there is a subsequence (not relabeled) such that
\begin{equation}\label{dod2}
P(R_\ep,Z_\ep)\to\overline{P(R,Z)}:=\overline{P(R,Z)}_n\ \mbox{weakly-* in} \ L^\infty(I\times\Omega),
\end{equation}
where (since $P$ is continuous and convexe, and since, in particular, $r_\ep\to r$ weakly in $L^1(I\times\Omega)$)
$$
P(R,Z):=P(R_n,Z_n)\le \overline{P(R,Z)}_n\ \mbox{a.e. in }\ I\times\Omega
$$
This is enough to pass to the limit in the momentum equation \eqref{E5}, in order to obtain:
\begin{equation} \label{E20}
\begin{split}
\left[ \intO{ (R+Z)\vu \cdot \bfphi } \right]_{t=0}^{t = \tau} &=
\int_0^\tau \intO{ \Big[ (R+Z) \vu \cdot \partial_t \bfphi + (R+Z) \vu \otimes \vu : \Grad \bfphi \\
&
+ \Ov{P(R,Z)} \Div \bfphi - \mathbb{S}(\Grad\vu)  : \Grad \bfphi \Big] }
\end{split}
\end{equation}
for any $\bfphi \in C^1([0,T]; X_n)$.

\subsubsection{Limit in the energy inequality}

We shall treat $H$ similarly as  $P$: as in (\ref{dod2})

$$
H(R_{\ep}, Z_{\ep})\to\ \mbox{weakly-* in}\ \overline{H(R,Z)}:=\overline{H(R, Z)}_n  \mbox{ in $L^\infty(I\times\Omega)$}.
$$
Since $H$ is convex continuous, since, in particular $r_\ep\to r$   weakly in $L^1(I\times\Omega)$, we have
\begin{equation}\label{dod1}
 \ 0\le H(R_n,Z_n)\le\overline{H(R,Z)}_n,
\end{equation}
where
$$
\|\overline{H(R, Z)}_n\|_{L^\infty(I;L^1(\Omega))},\;\mbox{is bounded uniformly with $n$}.
$$

Seeing (\ref{E15++}) and convexity of $H$ we may conclude that
$$
\int_0^\tau \int_{\Gamma^{\rm out}}H(R,Z)\vu_B\cdot\vc n{\rm d}S_x{\rm d}t\le\liminf_{\ep\to 0} \int_0^\tau \int_{\Gamma^{\rm out}}H(R_\ep,Z_\ep)\vu_B\cdot\vc n{\rm d}S_x{\rm d}t.
$$

Finally recalling (\ref{dod0}), we can pass to the limit also in the energy balance (\ref{E7}), where we obtain after omitting at the left hand side several non negative terms, the following energy inequality:

\begin{equation} \label{E7+}
\begin{split}
&\left[ \intO{\left[ \frac{1}{2} (R+Z)|\vv|^2 + \overline{ H(R,Z)} \right] } \right]_{t = 0}^{ t = \tau} +
\int_0^\tau \intO{\mathbb{S}(\Grad\vu):\Grad\vu } \dt
\\
&+\int_0^\tau \int_{\Gamma^{\rm out}} { H}(R,Z)  \vu_B \cdot \vc{n} \ {\rm d} S_x \dt
\leq
-
\int_0^\tau \intO{ \left[ (R+Z) \vu \otimes \vu + \overline{P(R,Z)} \mathbb{I} \right]  :  \Grad \vu_B } \dt \\
&+ \int_0^\tau { \intO{ (R+Z) \vu  \cdot\Grad \vu_B \cdot  \vu_B  } }
\dt
+ \int_0^\tau \intO{ \mathbb{S}(\Grad\vu) : \Grad \vu_B } \dt
\\
&- \int_0^\tau \int_{\Gamma^{\rm in}} { H}(R_B,Z_B)  \vu_B \cdot \vc{n} \ {\rm d} S_x \dt
\end{split}
\end{equation}
for a.a. $\tau\in I$

Before summarizing the results of this section we introduce the kinetic energy function which will be convenient to use in the last limit process.

For any $0\neq \xi\in\R^3$, we introduce the following convex lower semicontinuous function
\begin{equation}
\mathbb{E}_\xi[r, \vc m]: \R\times\R^d\ni[r, \vc m] \mapsto \left\{\begin{array}{c}
 \frac{|\vm \cdot \xi|^2}{r}\;\mbox{if $r>0$},
 \\
 0\;\mbox{if $r=0$, $\vm=0$},\\
 \infty\;\mbox{otherwise},
                   \end{array}\right.
\end{equation}
together with function
$$
\mathbb E [r, \vc m]: \R\times\R^d\ni[r, \vc m] \mapsto \left\{\begin{array}{c}
 \frac{\vm \otimes \vm}{r}\;\mbox{if $r>0$},
 \\
 0\;\mbox{if $r=0$, $\vm=0$},\\
 \infty\;\mbox{otherwise}.
                   \end{array}\right.
$$
and verify that with $R=R_n$, $Z=Z_n$, $\vc u=\vc u_n$, we have
$$
\mathbb E_\xi[R+Z, (R+Z)\vc u]=(R+Z)|\vu\cdot\xi|^2,\; \mathbb E_\xi[R+Z, (R+Z)\vc v]=(R+Z)|\vc v\cdot\xi|^2,\;\mathbb E_\xi(r,\vc m)=
\xi^T\mathbb{E}[r,\vc m]\xi.
$$
Finally, we denote
$$
\mathbb E_0=\mathbb E_{{\vc e}_1}+\ldots +\mathbb E_{{\vc e}_d},\;\mbox{ where ${\vc e}_i$
is canonical basis of $\R^d$}.
$$

\subsubsection{Conclusion for the limit $\ep\to 0$ }

To conclude, we summarize the result obtained in the limit $\ep\to 0$.

\begin{Proposition}[{{ Second level of approximate solutions ($n$ fixed)}}] \label{EP2}
Let $\Omega \subset \R^3$ be a bounded Lipschitz domain.
Let the data $(R_B,Z_B,\vu_B)$, $(R_0,Z_0,\vu_0)$  belong to the class { (\ref{Bdata2}), (\ref{Bbis})},
(\ref{Idata2}), (\ref{Ibis}). Suppose that assumptions (\ref{Pressure1}--\ref{Pressure2}) are satisfied.

Then for each fixed $n > 0$,  there exists
$(R_n,Z_n,\vu_n=\vv_n+\vu_B)$
in the class
$$
R,Z\in L^\infty( I\times\Omega),\
\partial_t(R,Z)\in L^2(I\times\Omega),\ R,Z\in L^\infty(I\times\partial\Omega),
$$
$$
\forall t\in\overline I,\
R(t),Z(t)> 0\ \mbox{a.e. in $\Omega$},\;
\mbox{for a.e. $t\in I$,}\
R(t),Z(t)> 0\ \mbox{a.e. in $\partial\Omega$},
$$
$$
\vv=C(\overline I; X_n), \;\partial_t\vv\in L^2(I;X_n)
$$
such that  the following holds:
\begin{enumerate}
\item Domination inequalities:
\begin{equation}\label{eq5.1+}
\begin{array}{c}
\forall t\in\overline I,\;
0<\underline c(n)\le
R_n(t,x),Z_n(t,x)\le {\overline c(n)} ,\; \underline b R_n(t,x) \leq Z_n(t,x) \leq \overline b R_n(t,x),
\mbox{a.e. in $\Omega$},
\\
\mbox{for a.a. $t\in I$},\;
0<\underline c(n)\le R_n(t,x), Z_n(t,x){ \le \overline c(n)},\; \underline b R_n(t,x) \leq Z_n(t,x) \leq \overline b R_n(t,x),
\mbox{a.e. in $\partial\Omega$};
\end{array}
\end{equation}
\item Continuity equations:
\begin{equation} \label{E22}
\begin{split}
&\left[ \intO{ r_n \varphi } \right]_{t = 0}^{t = \tau}=
\int_0^\tau \intO{ \Big[ r_n \partial_t \varphi + r_n \vu_n \cdot \Grad \varphi  \Big] }
\dt
- \int_0^\tau \int_{\Gamma^{\rm out}} \varphi r_n \vu_B \cdot \vc{n} \ {\rm d} S_x \dt \\
-&
\int_0^\tau \int_{\Gamma^{\rm in}} \varphi r_B \vu_B \cdot \vc{n}  \ {\rm d}S_x \dt,\
r(0, \cdot) = r_0,\; r\ \mbox{stands for}\ R,Z,
\end{split}
\end{equation}
for any $\varphi \in C^1([0,T] \times \Ov{\Omega})$;
\item Momentum equation:
\begin{equation} \label{E23}
\begin{split}
&\left[ \intO{ (R_n+Z_n) \vu \cdot \bfphi } \right]_{t=0}^{t = \tau} =
\int_0^\tau \int_\Omega \Big[ (R_n+Z_n) \vu_n \cdot \partial_t \bfphi \\
&+ \Big(\mathbb{E}[R_n+Z_n, \vc m_n] + \Ov{P(R,Z)}_n\mathbb I\Big) : \Grad \bfphi \\
& - \mathbb{S}(\Grad\vu_n)  : \Grad \bfphi \Big] {\rm d}x{\rm d}t,\ (R+Z)\vu(0)=\vc m_0,\;\vc m_n=(R_n+Z_n)\vu_n
\end{split}
\end{equation}
for any $\bfphi \in C^1([0,T]; X_n)$, where
\begin{equation}\label{E23+}
0\le P(R_n,Z_n)\le \Ov{P(R,Z)}_n,\ \|\Ov{P(R,Z)}_n\|_{L^\infty(I;L^1(\Omega))}\ \mbox{uniformly bounded}.
\end{equation}
\item The approximate energy inequality
\begin{equation} \label{E24}
\begin{split}
&\intO{\left[ \frac{1}{2} \mathbb{E}_0[R_n(\tau)+Z_n(\tau),\vc q_n(\tau)] + \overline{{ H}(R,Z)}_n  \right] }\\
&{-}\intO{\left[ \frac{1}{2} \mathbb{E}_0[R_0+Z_0(\tau),\vc q_0] + { H}(R_0,Z_0)  \right] }
\\
&+
\int_0^\tau \intO{\mathbb{S}(\Grad\vu_n):\Grad\vu_n } \dt
+\int_0^\tau \int_{\Gamma^{\rm out}} { H}(R_n,Z_n)  \vu_B \cdot \vc{n} \ {\rm d} S_x \dt
\ {\rm d}S_x \ \dt \\
&\leq
-
\int_0^\tau \intO{ \left[ \mathbb{E}[R_n+Z_n, \vc m_n] + \overline{P(R,Z)}_n \mathbb{I} \right]  :  \Grad \vu_B } \dt \\
&+ \int_0^\tau { \intO{ (R_n+Z_n) \vu_n  \cdot\Grad \vu_B \cdot  \vu_B  } }
\dt
+ \int_0^\tau \intO{ \mathbb{S}(\Grad\vu_n) : \Grad \vu_B } \dt
\\
&- \int_0^\tau \int_{\Gamma^{\rm in}} { H}(R_B,Z_B)  \vu_B \cdot \vc{n} \ {\rm d} S_x \dt,\;\vc q_n=(R_n+Z_n)\vc v_n,\;\vc q_0=
(R_0+Z_0)\vc v_0
\end{split}
\end{equation}
holds for a.a. $\tau\in I$, where $\vc q_n=\vv_n+\vu_B$.
\end{enumerate}

\end{Proposition}

\subsection{Existence for the bi-fluid system (limit $n\to\infty$).}

Our ultimate goal is to perform the limit $n \to \infty$ in the family of approximate solutions obtained in
Proposition \ref{EP2}.

\subsubsection{Limit in the continuity equation}\label{SCE}

 An easy application of Gronwall's lemma shows that the total energy represented by the expression on the left--hand side of the energy inequality \eqref{E24} remains bounded uniformly for $n \to \infty$. This, together  with the domination inequlities (\ref{eq5.1+}) yields estimates
 $$
 \|r_n\|_{L^\infty I;L^\gamma(\Omega))}\le c,\;\|r_n\|_{L^\gamma(I;
 L^\gamma(\Gamma^{\rm out};|\vu_B\cdot\vc n|{\rm d}S_x))}\le c
 $$
 $$
 \|r_n\vu_n\|_{L^\infty (I;L^{\frac{2\gamma}{\gamma+1}}(\Omega))}\le c,\;\|\Grad\vu_n\|_{L^2(I\times\Omega)}\le c.
 $$

 Consequently, extracting suitable subsequences if necessary, we may suppose
\begin{equation}\label{E29-}
r_n \to r \ \mbox{in}\ C_{{\rm weak}}([0,T]; L^\gamma(\Omega)),\
r_n\to r \ \mbox{weakly-(*) in}\ L^\infty(0,T; L^\gamma (\Gamma^{\rm out}; |\vu_B \cdot \vc{n}| \D S_x)),
\end{equation}
\begin{equation}\label{E29+1}
\vm_n = r_n \vu_n \to \vm \ \mbox{weakly-(*) in}\ L^\infty(0,T; L^{\frac{2 \gamma}{\gamma + 1}}(\Omega; \R^d)).
\end{equation}
and
\begin{equation}\label{E29+2}
\vu_n \to \vu \ \mbox{weakly in}\ L^2(0,T; W^{1,2}(\Omega; \R^d)).
\end{equation}
We show now that
\begin{equation} \label{E29}
\vm = r \vu \ \mbox{a.a. in}\ (0,T) \times \Omega.
\end{equation}
This is a direct application of (nontrivial) Lemma \ref{WANL1} (see Appendix), where we take
\[
r_n = r_n,\ v_n = u^i_n,\ i=1,\dots,d, \ \vc{g}_n = - r_n \vu_n,\ h_n = 0
\]
The hypotheses of the lemma are satisfied with exponents
$p = \gamma$, $r = s =
\frac{2 \gamma}{\gamma + 1}$

At this stage, we are able to perform the limit in the equations of continuity \eqref{E22} to obtain \eqref{W1}.

\subsubsection{Limit in the momentum equation}

The next step is to perform the same limit in the momentum equation \eqref{E23}.

Seeing, on one hand that $L^\infty(I;L^1(\Omega))\hookrightarrow
L^\infty(I; {\cal M}(\overline\Omega))$ and that
$(L^1(I;C(\overline\Omega))^*=L^\infty_{weak-*}(I,$ ${\cal M}(\overline\Omega))$, cf. Lemma \ref{Lemma3}, and on the other hand
that $P$ is convex, continuous and that (\ref{E29-}) holds
we get from (\ref{E23+}), in particular
$$
L^\infty_{weak-*}(I,{\cal M}^+(\overline\Omega))\ni {\mathfrak{R}}_1 :=\overline{\overline{P(R,Z)}}- P(R,Z),
$$
where $\overline{\overline{P(R,Z)}}$ is *-weak limit of a chosen
subsequence
$\overline{P(R,Z)}_n$ (not relabeled) in $L^\infty_{weak-*}(I;{\cal M}(\overline\Omega))$ (whose existence is guranteed by the Banach-Alaoglu-Bourbaki theorem).
Likewise,
$$
\mathbb{E}[R_n+Z_n,\vc m_n]\to\overline{\mathbb{E}[R+Z,\vc m]}\
\mbox{weakly-* in}\, L^\infty_{weak-*}(I,{\cal M}(\overline\Omega))
$$
where
$$
\mathfrak{R}_2=\overline{\mathbb{E}[R+Z,\vc m]}-\mathbb{E}[R+Z,\vc m]
\in L^\infty_{weak-*}(I,{\cal M}^+(\overline\Omega;\R^{d\times d}_{sym}))
$$
Indeed, since for all $0\neq \xi$, $\mathbb E_\xi$ is lower semicontinuous and convex, we have
$$
\xi^T\mathbb{E}[R+Z,\vc m]\xi\le \overline{\xi^T\mathbb{E}[R+Z,\vc m]\xi}, \;\xi\in \R^d.
$$
Thus, letting $n\to\infty$ in (\ref{E23}) we obtain the formulation (\ref{W3}) of the momentum equation
with
$$
\mathfrak{R}=\mathfrak{R}_1+\mathfrak{R}_2=
\overline{\mathbb{E}[R+Z,\vc m]}
-\mathbb{E}[R+Z,\vc m]+\Big(\overline{\overline{P(R,Z)}}- P(R,Z)\Big)\mathbb I.
$$

\subsubsection{Limit in the energy inequality}
By the same token, due to (\ref{dod1}),
$$
L^\infty_{weak-*}(I,{\cal M}^+(\overline\Omega))\ni {\mathfrak{E}}_1:=\overline{\overline{H(R,Z)}}- H(R,Z),
$$
where $\overline{\overline{H(R,Z)}}$ is *-weak limit of a chosen
subsequence
$\overline{H(R,Z)}_n$ (not relabeled) in $L^\infty_{weak-*}(I,{\cal M}(\overline\Omega))$.
Further, it is easy to see that
$$
\mathfrak{E}_2:=\overline{\mathbb E_0(R+Z,\vc q)}-{\mathbb E_0[R+Z,\vc q]}= \overline{\mathbb E_0(R+Z,\vc m)}-{\mathbb E_0[R+Z,\vc m]}
\in L^\infty_{weak-*}(I,{\cal M}^+(\overline\Omega)).
$$
Thus setting
$$
\mathfrak{E} (\tau)=\lim_{h\to 0+}\frac 1{2h}\int_{\tau-h}^{\tau+h}
\int_{\overline\Omega} {\rm d}(\mathfrak{E}_1(t)+\mathfrak{E}_2(t))
$$
(due to the Theorem on Lebesgue points, this limit is equal
to $\mathfrak{E}_1(\tau)+\mathfrak{E}_2(\tau)$ for a.a. $\tau\in I$. Due to (\ref{E29-}) and since $H$ is convex,
$$
\int_0^\tau \int_{\Gamma^{\rm out}}H(R,Z)\vu_B\cdot\vc n{\rm d}S_x{\rm d}t\le\liminf_{n\to \infty} \int_0^\tau \int_{\Gamma^{\rm out}}H(R_n,Z_n)\vu_B\cdot\vc n{\rm d}S_x{\rm d}t.
$$
We are ready to pass to the limit in energy inequality (\ref{E24})
in order to get (\ref{W5}).

\subsubsection{Compatibility conditions}\label{SCC}

By structural assumptions (\ref{Pressure1})--(\ref{Pressure2}) and convexity
$$
\overline a\Big(\overline{\overline{P(R,Z)}}-P(R,Z)\Big)\ge \overline{\overline{H(R,Z)}}-H(R,Z)
$$
and
$$
\underline a\Big(\overline{\overline{P(R,Z)}}-P(R,Z)\Big)\le \overline{\overline{H(R,Z)}}-H(R,Z)
$$
we deduce
$$
\min\{\frac 1 d,\frac 1{\overline a}\}\mathfrak{E}\le
\frac 1 d{\rm Tr} \mathfrak{R}\le \max\{\frac 1 d,\frac 1{\underline a}\}\mathfrak{E}.
$$
This implies compatibility conditions (\ref{W6}).

\subsubsection{Finite energy initial data}

At this stage we have proved existence of dissipative turbulent solutions
with the regular initial and boundary data in class \eqref{Bbis}--(\ref{Ibis})
and (\ref{Bdata2}), (\ref{Idata2}). In order to get finite energy initial data
\eqref{Idata1} and boundary data \eqref{Bdata1} (and (\ref{Bdata2}), (\ref{Idata2})) we have to perform the steps from Sections \ref{SCE}--\ref{SCC}
with boundary data $(r_{0,n},\vu_{0,n})$, and $(r_{B,n},\vu_B)$ in class
\eqref{Bbis}--(\ref{Ibis}), where $(r_{0,n},\vu_{0,n})$ and $r_{B,n}$
are approximations of $(r_{0},\vu_{0})$, and $r_{B}$ such that
\[
r_{0,n} \to r_0 \ \mbox{in}\ L^\gamma(\Omega),\; H(R_{0,n},Z_{0,n})\to
H(R_0,Z_0)\ \;\mbox{in}\ L^1(\Omega), r_{B,n}\to r_B\ \mbox{in}\ C(\overline\Omega),
\]
$$
(R_{0,n}+Z_{0,n})\vu_{0,n}\to \vc m_0\ \mbox{in}\ L^{\frac{2\gamma}{\gamma+1}}(\Omega),\ \mathbb{E}_0[R_{0,n}+Z_{0,n}, (R_{0,n}+Z_{0,n})\vu_{0,n}]\to
\mathbb{E}_0[R_0+Z_0,\vc m_0]\ \mbox{in}\  L^1(\Omega).
$$
In this way, we obtain
  \eqref{W2}, \eqref{W3}, \eqref{W5} with the desired finite energy initial data.

\section{Compatibility with classical solution: Proof of Theorem
\ref{ET2}}
\label{C}

In this section, we show Theorem \ref{ET2}: if a dissipative solution enjoys certain regularity, specifically if
\[
\vu \in C^1([0,T] \times \Ov{\Omega}; R^d),\ R,Z \in C^1([0,T] \times \Ov{\Omega}),\
\inf_{(0,T) \times \Omega} Z > 0,
\]
then $[R,Z,\vu]$ is a classical solution, meaning $\mathfrak{E} = \mathfrak{R} = 0$.

To see this, we realize that $(\vu - \vu_B)$ can be used as a test function in the momentum equation
\eqref{W3}, which, together with the equation of continuity \eqref{W1}, yield the total energy \emph{equality}:
\begin{equation} \label{C1}
\begin{split}
&\left[ \intO{\left[ \frac{1}{2} r |\vu - \vu_B|^2 + H(R,Z) \right] } \right]_{t = 0}^{ t = \tau}
+
\int_0^\tau \intO{ \mathbb{S}(\Grad\vu): \Grad \vu } \dt
\\
&{ +\int_0^\tau \int_{\Gamma^{\rm out}} H(R,Z)  \vu_B \cdot \vc{n} \ \D S_x \dt+ \int_0^\tau \int_{\Gamma^{\rm in}} H(R_B,Z_B)  \vu_B \cdot \vc{n} \ \D S_x \dt \	
}\\
&=
-
\int_0^\tau \intO{ \left[ (R+Z)\vu \otimes \vu + P(R,Z) \mathbb{I} \right]  :  \Grad \vu_B } \dt + \int_0^\tau \intO{ {(R+Z)} \vu  \cdot \Grad\vu_B \cdot \vu_B  }
\dt  \\ &+ \int_0^\tau \intO{ \mathbb{S} : \Grad \vu_B } \dt-
\int_0^\tau \int_{\Ov{\Omega}} \Grad (\vu_B - \vu) : \D \ \mathfrak{R}(t) \dt.
\end{split}
\end{equation}
Relation \eqref{C1} subtracted from the energy inequality \eqref{W5} gives rise to
\[
\int_{\Ov{\Omega}} 1 \D \ \mathfrak{E} (\tau)
\leq \int_0^\tau \int_{\Ov{\Omega}} \Grad \vu : \D \ \mathfrak{R}(t) \dt.
\]
Since for any fixed $i,j$, $\mathfrak{R}_{ij}$ is a signed Radon measure, we decompose $\mathfrak{R}_{ij}=\mathfrak{R}_{i,j}^++
\mathfrak{R}_{i,j}^-$; in particular $\pm\int_{\overline\Omega}{\rm d}\mathfrak{R}_{i,j}^\pm\ge 0$. Moreover, since $\mathfrak{R}\in {\cal M}^+
(\overline\Omega;R^{d\times d}_{\rm sym})$, we have for any $i$,
$\mathfrak{R}_{ii}\ge 0$; whence, in particular ${\rm Tr}(\mathfrak{R})
= {\rm Tr}(\mathfrak{R}^+)$ and ${\rm Tr}(\mathfrak{R}^-)=0$.
Consequently,
$$
\int_{\overline\Omega}\Grad\vu:{\rm d}\mathfrak{R}=\int_{\overline\Omega}\Grad\vu:{\rm d}\mathfrak{R}^+ -
\Grad\vu:{\rm d}(-\mathfrak{R}^-)\le
$$
$$
\sup_{I\times\Omega}|\Grad\vu| \mathbb{I}: \int_{\overline\Omega}({\rm d}\mathfrak{R}^+
+ {\rm d}(-\mathfrak{R}^-))
= \sup_{I\times\Omega}|\Grad\vu|\int_{\overline\Omega}{\rm d}{\rm Tr}\mathfrak{R}.
$$
This, together with the compatibility hypothesis \eqref{W6} and Gronwall lemma, yields the desired conclusion
$\mathfrak{E} = \mathfrak{R} = 0$.

\section{Relative energy: Proof of Theorem \ref{ET3}}
\label{RR}

The relative energy
$
\mathcal{E} \left(R, Z, \vu \ \Big|\ \mathfrak{r},\mathfrak{z}, \mathfrak{u} \right)$ can be rewritten as
\[
\begin{split}
&\frac{1}{2} (R+Z) |\vu - \mathfrak{u}|^2 + H(R,Z) -
\partial_RH(\mathfrak{r},\mathfrak{z}) (R - \mathfrak{r}) -
\partial_ZH(\mathfrak{r},\mathfrak{z}) (Z - \mathfrak{z})
-
P(\mathfrak{r},\mathfrak{z}) \\
&= \frac{1}{2} (R+Z) |\vu - \vu_B - (\mathfrak{u} - \vu_B)|^2  
+ H(R,Z) -
\partial_RH(\mathfrak{r},\mathfrak{z}) (R - \mathfrak{r}) -
\partial_ZH(\mathfrak{r},\mathfrak{z}) (Z - \mathfrak{z})
-
P(\mathfrak{r},\mathfrak{z})
\\&= \left[ \frac{1}{2} (R+Z) |\vu - \vu_B|^2 + H(,R,Z) \right] - (R+Z) \vu \cdot (\mathfrak{u} - \vu_B)
\\
&
 + \left[ \frac{1}{2} \Big( |\mathfrak{u}|^2 - |\vu_B|^2 \Big) -  \partial_RH(\mathfrak{r},\mathfrak{z}) \right] R
 + \left[ \frac{1}{2} \Big( |\mathfrak{u}|^2 - |\vu_B|^2 \Big) -  \partial_ZH(\mathfrak{r},\mathfrak{z}) \right] Z
+P(\mathfrak{r},\mathfrak{z}){ .}
\end{split}
\]
Our goal is to evaluate the time evolution of
\[
\intO{ \mathcal{E} \left(R,Z, \vu \ \Big|\ \mathfrak{r},\mathfrak{r}, \mathfrak{u} \right) }
\]
where $[R,Z, \vu]$ is a dissipative solutions and $[\mathfrak{r},\mathfrak{z}, \mathfrak{u}]$ are test functions in the class \eqref{W6}.

\medskip
\noindent
{\bf Step 1:}
\\
In accordance with the energy inequality \eqref{W5}, we get
\begin{equation} \label{RR1}
\begin{split}
&\left[ \intO{\left[ \frac{1}{2} (R+Z) |\vu - \vu_B|^2 + H(R,Z) \right] } \right]_{t = 0}^{ t = \tau} +
\int_0^\tau \intO{ \mathbb{S}(\Grad\vu):\Grad\vu } \dt
\\
& { +\int_0^\tau \int_{\Gamma^{\rm out}} H(R,Z)  \vu_B \cdot \vc{n} \ \D S_x \dt +\int_0^\tau \int_{\Gamma^{\rm in}} H(R_B,Z_B)  \vu_B \cdot \vc{n} \ \D S_x \dt}
+ \int_{\Ov{\Omega}} 1 \D \ \mathfrak{E} (\tau) \\	
\leq
&-
\int_0^\tau \intO{ P(R,Z)  \Div \vu_B } \dt +
\int_0^\tau \intO{ {(R+Z)}\Big(\vu\cdot\Grad\vu_B \cdot\vu_B
-\vu\cdot\Grad\vu_B\cdot\vu\Big) }
\dt\\
&+ \int_0^\tau \intO{ \mathbb{S}(\Grad\vu) : \Grad \vu_B } \dt -
\int_0^\tau \int_{\Ov{\Omega}} \Grad \vu_B : \D \ \mathfrak{R}(t) \dt{ .}
\end{split}
\end{equation}

\noindent
{\bf Step 2:}
\\
Plugging $\bfphi = \mathfrak{u} - \vu_B$ in the momentum equation \eqref{W3}, we get
\begin{equation} \label{RR2}
\begin{split}
&\left[ \intO{ (R+Z) \vu \cdot (\mathfrak{u} - \vu_B) } \right]_{t=0}^{t = \tau} = \int_0^\tau \int_\Omega \Big[ (R+Z) \Big(\vu \cdot \partial_t \mathfrak{u}  +  \vu \cdot\Grad\mathfrak{u}\cdot \vu -\vu\cdot \Grad\vu_B\cdot\vu\Big)
\\
&
+ P(R,Z)\Div(\mathfrak{u} - \vu_B) - \mathbb{S}(\Grad\vu) : \Grad (\mathfrak{u} - \vu_B) \Big] {\rm d}x \dt
+ \int_0^\tau \int_{\Ov{\Omega}} \Grad (\mathfrak{u} - \vu_B) : \D \ \mathfrak{R}(t) \dt{.}
\end{split}
\end{equation}
\medskip
\noindent
{\bf Step 3:}
\\
Finally, we consider $\varphi = \left[ \frac{1}{2} \Big( |\mathfrak{u}|^2 - |\vu_B|^2 \Big) -  \partial_RH(\mathfrak{r},\mathfrak{z}) \right]$ in
the equation of continuity \eqref{W1} with $r=R$
and  $\varphi = \left[ \frac{1}{2} \Big( |\mathfrak{u}|^2 - |\vu_B|^2 \Big) -  \partial_ZH(\mathfrak{r},\mathfrak{z}) \right]$
in
the equation of continuity \eqref{W1} with $r=Z$
obtaining:
\begin{equation} \label{RR3}
\begin{split}
&\left[ \intO{ r \left[ \frac{1}{2} \Big( |\mathfrak{u}|^2 - |\vu_B|^2 \Big) -  \partial_rH(\mathfrak{r},\mathfrak{z})) \right] } \right]_{t = 0}^{t = \tau}  \\&{ -
\int_0^\tau \int_{\Gamma^{\rm out}}\partial_rH(\mathfrak{r},\mathfrak{z})r \vu_B \cdot \vc{n} \ \D \ S_x
-
\int_0^\tau \int_{\Gamma^{\rm in}}\partial_rH(\mathfrak{r}, \mathfrak{z})r_B \vu_B \cdot \vc{n} \ \D \ S_x
}
\\ &=
\int_0^\tau \intO{ \Big[ r \partial_t\left( \frac{1}{2} |\mathfrak{u}|^2 -  \partial_rH(\mathfrak{r}, \mathfrak{z}) \right)  +
r \vu \cdot \Grad \left( \frac{1}{2} \Big( |\mathfrak{u}|^2 - |\vu_B|^2 \Big) -  \partial_rH(\mathfrak{r}, \mathfrak{z}) \right) \Big] } \dt,
\end{split}
\end{equation}
where $r$ stands once for $R$ and once for $Z$.

\medskip
\noindent
{\bf Step 4:}\\
Summing up \eqref{RR1},\eqref{RR2}, subtracting (\ref{RR3})$_{r=R}$ and (\ref{RR3})$_{r=Z}$, { and using identity
$$
\Big[\intO{\Big(R\partial_RH(\mathfrak{r},\mathfrak{z})+
Z\partial_ZH(\mathfrak{r},\mathfrak{z})-H(\mathfrak{r},\mathfrak{z})\Big)}\Big]_0^{\tau}=\int_0^\tau\intO{\partial_tP(\mathfrak{r},\mathfrak{z})},
$$
}
we get
\begin{equation} \label{RR4}
\begin{split}
&\left[ \intO{\mathcal{E}\left( R,Z, \vu \ \Big|\ \mathfrak{r},
\mathfrak{z},\mathfrak{u} \right) } \right]_{t = 0}^{ t = \tau}
\\
&+
\int_0^\tau \intO{ \mathbb{S}(\Grad\vu) :\Grad\vu} \dt - \int_0^\tau \intO{ \mathbb{S}(\Grad\vu) : \Grad \mathfrak{u} } \dt \\
&+\int_0^\tau \int_{\Gamma^{\rm out}} \left[ H(R,Z) - R\partial_RH(\mathfrak{r},\mathfrak{z}) -Z\partial_ZH(\mathfrak{r},\mathfrak{z}) \right]  \vu_B \cdot \vc{n} \ \D S_x \dt
\\
&+\int_0^\tau \int_{\Gamma^{\rm in}} \left[ H(R_B,Z_B) - R_B\partial_RH(\mathfrak{r},\mathfrak{z}) -Z_B\partial_ZH(\mathfrak{r},\mathfrak{z}) \right]  \vu_B \cdot \vc{n} \ \D S_x \dt
\\
&+ \int_{\Ov{\Omega}} 1 \D \ \mathfrak{E} (\tau) \\	
\leq
&- \int_0^\tau \intO{ \Big[ (R+Z)( \vu \cdot \partial_t \mathfrak{u}  + \vu\cdot\Grad\mathfrak{u} \cdot\vu )
+ P(R,Z) \Div \mathfrak{u}  \Big] } \dt \\
&- \int_0^\tau \int_{\Ov{\Omega}} \Grad \mathfrak{u} : \D \ \mathfrak{R}(t) \ \dt  + \int_0^\tau \intO{ \partial_t P(\mathfrak{r},\mathfrak{z}) } \dt
\\
&+\int_0^\tau \intO{ \Big[ R \partial_t\left( \frac{1}{2} |\mathfrak{u}|^2 -  \partial_RH(\mathfrak{r},\mathfrak{z}) \right)  +
R \vu \cdot \Grad \left( \frac{1}{2} |\mathfrak{u}|^2  -  \partial_RH(\mathfrak{r},\mathfrak{z}) \right) \Big] } \dt
\\
&
+\int_0^\tau \intO{ \Big[ Z\partial_t\left( \frac{1}{2} |\mathfrak{u}|^2 -  \partial_ZH(\mathfrak{r},\mathfrak{z}) \right)  +
Z \vu \cdot \Grad \left( \frac{1}{2} |\mathfrak{u}|^2  -  \partial_ZH(\mathfrak{r},\mathfrak{z}) \right) \Big] } \dt{.}
\end{split}
\end{equation}
\medskip
\noindent
{\bf Step 5:}\\
Now we regroup conveniently the terms in (\ref{RR4}). To this end we proceed as follows:
\begin{enumerate}
 \item Since
 $$
 \partial_tP(\mathfrak{r},\mathfrak{z})=\partial_RP(\mathfrak{r},\mathfrak{z})\partial_t\mathfrak{r}+\partial_ZP(\mathfrak{r},\mathfrak{z})
 \partial_t\mathfrak{z} +\Div(P(\mathfrak{r},\mathfrak{z})\mathfrak{u})
 +(\mathfrak{r}\partial_RP(\mathfrak{r},\mathfrak{z})+ \mathfrak{z}\partial_ZP(\mathfrak{r},\mathfrak{z}))\Div\mathfrak{u}
 $$
 $$
 -\Div(P(\mathfrak{r},\mathfrak{z})\mathfrak{u})-(\mathfrak{z}\partial_ZP(\mathfrak{r},\mathfrak{z})+\mathfrak{z}\partial_ZP(\mathfrak{r},\mathfrak{z}))\Div\mathfrak{u}{,}
 $$
 {we} have
$$
\intO{\partial_tP(\mathfrak{r},\mathfrak{z})}=
\intO{\partial_RP(\mathfrak{r},\mathfrak{z})(\partial_t\mathfrak{r}+\Div(\mathfrak{r}\mathfrak{u}))}
+\intO{\partial_ZP(\mathfrak{r},\mathfrak{z})(\partial_t\mathfrak{z}+\Div(\mathfrak{z}\mathfrak{u}))}
$$
$$
-\int_{\partial\Omega}P(\mathfrak{r},\mathfrak{z})\mathfrak{u}\cdot\vc n{\rm d}S_x
-\intO{(\mathfrak{r}\partial_RP(\mathfrak{r},\mathfrak{z})+\mathfrak{z}\partial_ZP(\mathfrak{r},\mathfrak{z} ))\Div\mathfrak{u}}+\intO{P(\mathfrak{r},\mathfrak{z})\Div\mathfrak{u}},
$$
where
$$
-\int_{\partial\Omega}P(\mathfrak{r},\mathfrak{z})\mathfrak{u}\cdot\vc n{\rm d}S_x=
\int_{\partial\Omega}\Big(H(\mathfrak{r},\mathfrak{z})-\mathfrak    {r}\partial_RH((\mathfrak{r},\mathfrak{z})-
\mathfrak{z}\partial_ZH(\mathfrak{r},\mathfrak{z})\Big)
{\vu}_B\cdot\vc n{\rm d}S_x.
$$

\item We calculate
$$
R\Big(\partial_t\partial_RH(\mathfrak{r},\mathfrak{z})+\vu\cdot\Grad \partial_RH(\mathfrak{r},\mathfrak{z})\Big)
$$
$$
=R\partial_R^2H(\mathfrak{r},\mathfrak{z})\Big(\partial_t\mathfrak{r}+\Div(\mathfrak{r}\mathfrak{u})\Big)
-R\partial^2_R H(\mathfrak{r},\mathfrak{z})\mathfrak{r}\Div\mathfrak{u}
+R\partial^2_R H(\mathfrak{r},\mathfrak{z})(\vu-\mathfrak{u})\cdot\nabla\mathfrak{r}
$$
$$
+ R\partial_R\partial_ZH(\mathfrak{r},\mathfrak{z})\Big(\partial_t\mathfrak{r}+\Div(\mathfrak{r}\mathfrak{u})\Big)
-R\partial_R\partial_Z H(\mathfrak{r},\mathfrak{z})\mathfrak{z}\Div\mathfrak{u}
+R\partial_R\partial_Z H(\mathfrak{r},\mathfrak{z})(\vu-\mathfrak{u})\cdot\nabla\mathfrak{z}.
$$
Similar calculation holds for
$$
Z\Big(\partial_t\partial_ZH(\mathfrak{r},\mathfrak{z})+\vu\cdot\Grad \partial_ZH(\mathfrak{r},\mathfrak{z})\Big).
$$
Consequently,
$$
R\Big(\partial_t\partial_RH(\mathfrak{r},\mathfrak{z})+\vu\cdot\Grad \partial_RH(\mathfrak{r},\mathfrak{z})\Big)+
Z\Big(\partial_t\partial_ZH(\mathfrak{r},\mathfrak{z})+\vu\cdot\Grad \partial_ZH(\mathfrak{r},\mathfrak{z})\Big)
$$
$$
=\Big(R\partial_R^2H(\mathfrak{r},\mathfrak{z})+
Z\partial_R\partial_ZH(\mathfrak{r},\mathfrak{z})\Big)\Big(\partial_t\mathfrak{r}+\Div(\mathfrak{r}\mathfrak{u})\Big)
+\Big(Z\partial_Z^2H(\mathfrak{r},\mathfrak{z})+
R\partial_R\partial_ZH(\mathfrak{r},\mathfrak{z})\Big)\Big(\partial_t\mathfrak{z}+\Div(\mathfrak{z}\mathfrak{u})\Big)
$$
$$
-R\partial_RP(\mathfrak{r},\mathfrak{z})\Div\mathfrak{u}
-Z\partial_ZP(\mathfrak{r},\mathfrak{z})\Div\mathfrak{u}
$$
$$
+(R-\mathfrak{r})\partial^2_RH(\mathfrak{r},\mathfrak{z})(\vu-\mathfrak{u})\cdot\nabla\mathfrak{r}
+(R-\mathfrak{r})\partial_R\partial_Z H(\mathfrak{r},\mathfrak{z})(\vu-\mathfrak{u})\cdot\nabla\mathfrak{z}
$$
$$
+(Z-\mathfrak{z})\partial^2_ZH(\mathfrak{r},\mathfrak{z})(\vu-\mathfrak{u})\cdot\nabla\mathfrak{z}
+(Z-\mathfrak{z})\partial_R\partial_Z H(\mathfrak{r},\mathfrak{z})(\vu-\mathfrak{u})\cdot\nabla\mathfrak{r}
$$
$$
+\partial_R P(\mathfrak{r},\mathfrak{z})\nabla \mathfrak{r}\cdot(\vu-\mathfrak{u}) +\partial_ZP(\mathfrak{r},\mathfrak{z})\nabla \mathfrak{z}\cdot(\vu-\mathfrak{u}),
$$
where we have used the property (\ref{H+}) in several lines.
\item We also have
$$
(R+Z)(\mathfrak{u}-\vu)\cdot\partial_t\mathfrak{u}=
\frac{R+Z}{\mathfrak{r}+\mathfrak{z}}\partial_t((\mathfrak{r}+\mathfrak{z})\mathfrak{u})\cdot (\mathfrak{u}-\vu)
-\frac{R+Z}{\mathfrak{r}+\mathfrak{z}}
(\mathfrak{u}-\vu)\cdot\mathfrak{u}\Big(\partial_t(\mathfrak{r}+\mathfrak{z})+\Div(({\mathfrak{r}+\mathfrak{z}})\mathfrak{u})\Big)
$$
$$
+\frac{R+Z}{\mathfrak{r}+\mathfrak{z}}
(\mathfrak{u}-\vu)\cdot\mathfrak{u}\Div(({\mathfrak{r}+\mathfrak{z}})\mathfrak{u})
$$
and
$$
(R+Z)\vu\cdot\Grad\mathfrak{u}\cdot(\mathfrak{u}-\vu)=
\frac{R+Z}{\mathfrak{r}+\mathfrak{z}} (\mathfrak{r}+\mathfrak{z})\mathfrak{u}\cdot\Grad\mathfrak{u}\cdot(\mathfrak{u}-\vu) -(R+Z)(\mathfrak{u}-\vu)\cdot\Grad\mathfrak{u}\cdot (\mathfrak{u}-\vu);
$$
whence
$$
(R+Z)(\mathfrak{u}-\vu)\partial_t\mathfrak{u}+
(R+Z)\vu\cdot\Grad\mathfrak{u}\cdot(\mathfrak{u}-\vu)
$$
$$
=
\frac{R+Z}{\mathfrak{r}+\mathfrak{z}}\Big(
\partial_t((\mathfrak{r}+\mathfrak{z})\mathfrak{u})+
\Div((\mathfrak{r}+\mathfrak{z})\mathfrak{u}\otimes\mathfrak{u})\Big)\cdot(\mathfrak{u}-\vu)
$$
$$
-(R+Z)(\mathfrak{u}-\vu)\cdot\Grad\mathfrak{u}\cdot (\mathfrak{u}-\vu)
-\frac{R+Z}{\mathfrak{r}+\mathfrak{z}}
(\mathfrak{u}-\vu)\cdot\mathfrak{u}\Big(\partial_t(\mathfrak{r}+\mathfrak{z})+\Div(({\mathfrak{r}+\mathfrak{z}})\mathfrak{u})\Big){.}
$$
\end{enumerate}

This calculations guide us in regrouping the terms in (\ref{RR4}) in order to obtain (\ref{RR5}). Theorem \ref{ET3}
is thus proved.

\section{Weak--strong uniqueness: Proof of Theorem \ref{ET4}}
\label{WU}

Our goal is to show that any dissipative solution coincides with the strong solution emanating from the same initial data { and boundary conditions}.
Assuming the strong solution $[\mathfrak{r},\mathfrak{z}, \mathfrak{u}]$--which solves \eqref{P1}--\eqref{P7} in the classical sense-- belongs to the class \eqref{RR6}, the obvious idea is to use the relative energy
inequality \eqref{RR5}.  We may rewrite
\eqref{RR5} as
\begin{equation} \label{WU1}
\begin{split}
&\left[ \intO{\mathcal{E}\left( R, Z, \vu \ \Big|\ \mathfrak{r}, \mathfrak{z},\mathfrak{u} \right) } \right]_{t = 0}^{ t = \tau}+
\int_0^\tau \intO{ \mathbb{S}( \Grad (\vu-\mathfrak{u})):\Grad(\vu-\mathfrak{u})}\dt \\
&+\int_0^\tau \int_{\Gamma^{\rm out}} \left[ H(R,Z) - \partial_RH(\mathfrak{r},\mathfrak{z}) (R - \mathfrak{r}) -
\partial_ZH(\mathfrak{r},\mathfrak{z}) (Z - \mathfrak{r})-
H(\mathfrak{r},\mathfrak{z})  \right]  \vu_B \cdot \vc{n} \ \D S_x \dt
\\
&+ \int_{\Ov{\Omega}} 1 \D \ \mathfrak{E} (\tau)
\leq - \int_0^\tau \intO{ (R+Z) (\mathfrak{u} - \vu)  \cdot \Grad \mathfrak{u}\cdot (\mathfrak{u} - \vu) } \dt\\
&-
\int_0^\tau \intO{ \Big[ P(\mathfrak{r},\mathfrak{z}) - \partial_RP(\mathfrak{r},\mathfrak{z}) (R - \mathfrak{r}) -
\partial_ZP(\mathfrak{r},\mathfrak{z}) (Z - \mathfrak{r})
-P(\mathfrak{r},\mathfrak{z}) \Big] \Div \mathfrak{u} } \dt
\\
&+ \int_0^\tau \intO{ \Big(\frac{R+Z}{\mathfrak{r}+\mathfrak{z}} -1\Big)(\mathfrak{u} - \vu) \cdot \Big[ \partial_t ((\mathfrak{r}+\mathfrak{z}) \mathfrak{u})  +  \Div ((\mathfrak{r}+\mathfrak{z}) \mathfrak{u} \otimes \mathfrak{u})
   \Big] } \dt
  \\
  &+ \int_0^\tau \intO{
(R-\mathfrak{r})(\mathfrak{u}-\vu)\cdot\Big(\nabla\mathfrak{r}\partial^2_RH(\mathfrak{r},\mathfrak{z})
+\nabla\mathfrak{z}\partial_R\partial_Z H(\mathfrak{r},\mathfrak{z})\Big)}{\rm d}t
\\
&+\int_0^\tau \intO{
(Z-\mathfrak{z})(\mathfrak{u}-\vu)\cdot\Big(\nabla\mathfrak{z}\partial^2_ZH(\mathfrak{r},\mathfrak{z})
+\nabla\mathfrak{r}\partial_R\partial_Z H(\mathfrak{r},\mathfrak{z})\Big)
  } \dt
  \\
&- \int_0^\tau \int_{\Ov{\Omega}} \Grad \mathfrak{u} : \D \ \mathfrak{R}(t) \dt{ ,}
\end{split}
\end{equation}
where we have used equation \eqref{P2}
for $\mathfrak{r}$, $\mathfrak{z}$, $\mathfrak{u}$ in the seventh line of formula (\ref{RR5}) and
the identity
\[
\int_0^\tau \intO{  (\mathfrak{u} - \vu) \cdot \Div
{\mathbb{S}}(\Grad\mathfrak{u}) } \dt =
- \int_0^\tau \intO{ \mathbb{S}(\Grad\mathfrak{u}): \Grad(\mathfrak{u} - \vu)   } \dt.
\]

Let
$$
0<\underline r<\underline{\mathfrak{r}}:=\inf_{I\times\Omega}\mathfrak{ r}<\overline{\mathfrak{r}}:=\sup_{I\times\Omega}\mathfrak r<\overline r<\infty,\
0<\underline z<\underline{\mathfrak{z}}:=\inf_{I\times\Omega}\mathfrak{z}<\overline{\mathfrak{r}}:=\sup_{I\times\Omega}\mathfrak{ z}<\overline z<\infty,
$$
and denote
$$
K=[\underline{\mathfrak{r}},\overline{\mathfrak{r}}]\times [\underline{\mathfrak{z}},\overline{\mathfrak{z}}]\cap\overline{{\cal O}},\ L=(\underline r,\overline r)\times
(\underline z,\overline z)\cap {\cal O}
$$

It follows from the structural hypothesis \eqref{Pressure2} that
$E(R,Z|\mathfrak{r},\mathfrak{z})$ majorates $R+Z +1$ and $P(R,Z)$
outside $L$, and $(R-\mathfrak{ r})^2+
(Z-\mathfrak{z})^2$ in $\overline L$ for any $(\mathfrak{r},
\mathfrak{z})\in K$ - with a multiplicative constant dependent solely of $L$ and $K$. Adding to these ingredients the compatibility condition (\ref{W6}) we verify that the absolute value of the r.h.s. of (\ref{WU1})
is bounded by
\[
 c\left(\underline{\mathfrak{z}},\| \mathfrak{r}, \mathfrak{z}, \mathfrak{u} \|_{C^1(\overline I\times\overline\Omega)},\delta \right) \left[
\int_0^\tau \intO{ \mathcal{E}\left(R,Z, \vu \ \Big|\ \mathfrak{r},\mathfrak{z}, \mathfrak{u} \right) } \dt +
\int_0^\tau \left(\int_{\Ov{\Omega}} \D \ \mathfrak{E} (t) \right) \dt     \right]+\delta\|\vu-\mathfrak{u}\|_{L^2(I\times\Omega)}^2
\]
with any $\delta>0$.

The term $\delta\|\vu-\mathfrak{u}\|_{L^2(I\times\Omega)}^2$
can be ``absorbed'' at the left hand side by the term
$$
\int_0^\tau \intO{ \mathbb{S}( \Grad (\vu-\mathfrak{u})):\Grad(\vu-\mathfrak{u})}\dt\ge c\|\vu-\mathfrak{u}\|_{L^2(I;W^{1,2}(\Omega))}^2-
\int_0^\tau\intO{(R+Z)(\vu-\mathfrak{u})^2}{\rm d}t,
$$
where the latter is true by virtue of the Korn and Sobolev inequalities (cf.\cite[Theorem 10.17]{FeNoB}).

Inequality (\ref{WU1}) therefore becomes
\begin{equation}\label{WU5}
\intO{{\cal E}(R,Z,\vu|\mathfrak{r},\mathfrak{z},\mathfrak{u})(\tau,\cdot)}+\int_{\overline\Omega}
{\rm d}\mathfrak{E}(\tau)
\le
c \int_0^\tau\intO{{\cal E}(R,Z,\vu|\mathfrak{r},\mathfrak{z},\mathfrak{u})}{\rm d}t+
\int_0^\tau\int_{\overline\Omega}{\rm d}\mathfrak{E}(t){\rm d}t.
\end{equation}

Thus applying the standard Gronwall argument to \eqref{WU5} we obtain the desired conclusion:
\[
r = \mathfrak{r}, \ \vu = \mathfrak{u}, \ \mathfrak{R} = \mathfrak{E} = 0.
\]

\section{Appendix}

We recall here some key  lemmas which we have used in the proofs.
The first lemma is proved in \cite[Lemma 8.1]{AbFeNo}.

\begin{Lemma} \label{WANL1}

Let $Q = (0,T) \times \Omega$, where $\Omega \subset R^d$ is a bounded domain. Suppose that
\[
r_n \to r \ \mbox{weakly in}\ L^p(Q), \ v_n \to v \ \mbox{weakly in}\ L^q(Q),\ p > 1, q > 1,
\]
and
\[
r_n v_n \to w \ \mbox{weakly in}\ L^r(Q), \ r >  1.
\]
In addition, let
\[
\partial_t r_n = \Div \vc{g}_n + h_n\ \mbox{in}\ \mathcal{D}'(Q),\ \| \vc{g}_n \|_{L^s(Q; R^d)} \aleq 1,\ s > 1,\
h_n \ \mbox{precompact in}\ W^{-1,z},\ z > 1,
\]
and
\[
\left\| \Grad v_n \right\|_{\mathcal{M}(Q; R^d)} \aleq 1 \ \mbox{uniformly for}\ n \to \infty.
\]
Then
\[
w = r v \ \mbox{a.a. in}\ Q.
\]
\end{Lemma}

The second lemme in Lemma 2.11 and Corollary 2.2 in
Feireisl \cite{EF70}.

\begin{Lemma}\label{Lemma2}
Let $O \subset \R^d$, $d\ge 2$, be a measurable set and $\{ \vc{v}_n
\}_{n=1}^{\infty}$ a sequence of functions in $L^1(O; \R^M)$ such
that
$$
\vc{v}_n \rightharpoonup \vc{v} \ \mbox{ in}\ L^1(O; \R^M).
$$
Let $\Phi: R^M \to (-\infty, \infty]$ be a lower semi-continuous
convex function such that $\Phi(\vc v_n)$ is bounded in $L^1({ O})$.

Then $\Phi(\vc v):O\mapsto R$ is integrable and
$$
\int_{O} \Phi(\vc v){\rm d} x\le \liminf_{n\to\infty} \int_{O} \Phi(\vc v_n){\rm d} x.
$$
\end{Lemma}

{ The last lemma we wish to recall deals with the duals of Bochner spaces.}
Let $X$  be a Banach space. For $1\le p<\infty$, we introduce Bochner-type spaces
$$
L^p(I,X)=\mbox{$\{f:I\to X$ measurable, $\|f\|^p_{L^p(I;X)}:=\int_I\|f(t)\|^p_X{\rm d}t<\infty\}$},
$$
$$
L_{\rm weak-*}^p(I,X^*)=
\mbox{$\{f:I\to X^*$ weakly-* measurable,}
$$
$$
\mbox{$\|f(\cdot)\|_{X^*}$ measurable,
$\|f\|^p_{L^p(I;X^*)}:=\int_I\|f(t)\|^p_{X^*}{\rm d}t<\infty\}$}.
$$
They are Banach spaces with corresponding norms
$\|f\|^p_{L^p(I;X)}$ resp; $\|f\|^p_{L^p(I;X^*)}$.
It is not true in general that $(L^p(I,X))^*$ can be identified
with $L^{p'}(I,X^*)$, $1/p+1/p'=1$. However, the following lemma holds,
cf. Pedregal\cite{Pedregal}:
\begin{Lemma}\label{Lemma3}
 Let $X$ be a separable Bancah space. Then
 $$(L^p(I,X))^*=
L^{p'}_{\rm weak-*}(I,X^*)$$
under the duality mapping
$$
<f,g>=\int_I<f(t),g(t)>_{X^*,X}{\rm d}t.
$$
\end{Lemma}

The particular case we are interested in in this paper
deal with
$$
X=C(\overline\Omega),\;{\mbox{whence}}\;
X^*={\cal M}(\overline\Omega),
$$
cf. Rudin\cite[Theorem 2.14]{Rudin}.

\def\cprime{$'$} \def\ocirc#1{\ifmmode\setbox0=\hbox{$#1$}\dimen0=\ht0
  \advance\dimen0 by1pt\rlap{\hbox to\wd0{\hss\raise\dimen0
  \hbox{\hskip.2em$\scriptscriptstyle\circ$}\hss}}#1\else {\accent"17 #1}\fi}


\end{document}